\documentclass{IJFS}

\usepackage[utf8]{inputenc}
\usepackage{amsfonts,amsthm,amsmath,amssymb,dsfont}
\usepackage{amsmath,    enumitem}
\usepackage{graphicx}
\usepackage[all]{xy}
\usepackage{subfigure}
\usepackage{tikz}
\usetikzlibrary[arrows,shapes,graphs,decorations.pathmorphing,backgrounds,positioning,fit,petri]
\usepackage{amscd, eufrak}

 \usepackage{amssymb}

 \usepackage{mathrsfs}
 \usepackage{amsfonts}

 \usepackage{verbatim}

\usepackage{color,colortbl}

\definecolor{OrcidGreen}{rgb}{0.6,0.8,0.2}

 \usepackage[singlelinecheck=false]{caption}
 \usepackage[linesnumbered,ruled]{algorithm2e}
 \usepackage[colorlinks=true,linkcolor=black,citecolor=black,urlcolor=blue]{hyperref}

\usepackage[square,numbers,sort&compress]{natbib}
\usepackage{url}
\usepackage{hyphenat}
\usepackage{parskip}

\usepackage{tabularray}

\usepackage{upgreek}
\usepackage{breqn}
\usepackage{mathtools}

\usepackage{pifont}
\usepackage{fancyhdr}
\usepackage{subfiles}
\usepackage{fontawesome5}

\usepackage{amsthm}

\newtheorem{theorem}{Theorem}[section]
\newtheorem{lemma}[theorem]{Lemma}

\newtheorem{corollary}[theorem]{Corollary}

\newtheorem{definition}[theorem]{Definition}
\newtheorem{example}[theorem]{Example}
\newtheorem{remark}[theorem]{Remark}

\newtheorem{step}[theorem]{Step}

\begin{document}

\addtocounter{page}{0}

\titleijfs{On the resolution and linear optimization problems subject to a system of bipolar fuzzy relational equalities defined with continuous Archimedean t-norms}
\author[1 \href{https://orcid.org/0000-0002-9224-8470}{\textcolor{OrcidGreen}{\faOrcid}}]{A. Ghodousian}
\author[2 \href{https://orcid.org/0009-0006-4046-2674}{\textcolor{OrcidGreen}{\faOrcid}}]{M. Sedigh Chopannavaz}
\author[3,4,5 \href{https://orcid.org/0000-0002-9335-9930}{\textcolor{OrcidGreen}{\faOrcid}}]{W. Pedrycz}
\affil[1]{Faculty of Engineering Science, College of Engineering, University of Tehran, P.O.Box 11365-4563, Tehran, Iran.}
\affil[2]{Department of Engineering Science, College of Engineering, University of Tehran, Tehran, Iran.}
\affil[3]{Department of Electrical and Computer Engineering, University of Alberta, Edmonton, AB T6G 2R3, Canada}
\affil[4]{Institute of Systems Engineering, Macau University of Science and Technology, Taipa 999078, Macau SAR, China}
\affil[5]{Research Center of Performance and Productivity Analysis, Istinye University, Istanbul, Türkiye}

\emails{a.ghodousian@ut.ac.ir, chopannavaz@ut.ac.ir and wpedrycz@ualberta.ca}
\CorrespondAuthor{A. Ghodousian}
\oddPageHead{On the resolution and linear optimization problems under BFRE system defined with continuous Archimedean t-norms}
\evenPageHead{A. Ghodousian, M. S. Chopannavaz and W. Pedrycz}

\abstractijfs{
Bipolar fuzzy relational equations are an interesting generalization of fuzzy relational equations that occur prominently in information processing, possibility theory, and preference modeling. In the latest related research, the authors studied max-strict bipolar fuzzy relational equations by categorizing the constraints of the problem into four groups, each requiring a different approach to determine the corresponding feasible solution set, and finding a global optimum by selecting a point having the minimum objective value among all the local optimal solutions. This paper considers the linear objective function optimization with respect to a more general class of bipolar fuzzy relational equations, where the fuzzy compositions are defined by an arbitrary continuous Archimedean t-norm. In addition, a faster method for finding a global optimum is proposed that, unlike the previous work, does not require obtaining all local optimal solutions and classifying the constraints (and therefore, it does not require applying different approaches to check the feasibility of constraints in different groups and the optimality of solutions in their feasible regions). Analytical concepts and properties of the Archimedean bipolar fuzzy equations are investigated and two necessary conditions are presented to conceptualize the feasibility of the problem. It is shown that the feasible solution set can be resulted by a union of the finite number of compact sets, where each compact set is obtained by a function (called admissible function in this paper). Moreover, to accelerate identification of the mentioned compact sets (and therefore, to speed up solution finding), four simplification techniques are presented, which are based on either omitting redundant constraints and/or eliminating unknowns by assigning them a fixed value. Also, three additional simplification techniques are given to reduce the search domain by removing some parts of the feasible region that do not contain optimal solutions. Subsequently, a method is proposed to find an optimal solution for the current linear optimization problems. The proposed method consists of two accelerative strategies that are used during the problem solving process. By the first strategy, the method neglects some candidate solutions that are not optimal, by considering only a subset of admissible functions (called modified functions in this paper). As for the second strategy, a branch-and-bound method is used to delete non-optimal branches. Then, the method is summarized in an algorithm that represents all essential steps of the solution and finally, the whole method is applied in an example that has been chosen in such a way that the various situations are illustrated.
}

\keywordsijfs{Bipolar fuzzy relational equations, Archimedean t-norms, Max-Archimedean compositions, Global optimization, Linear optimization.
}

\Vol{22}	\No{4}	\Year{2025}	\Pages{137-160}
\Received{January 2025}	\Revised{July 2025}		\Accepted{August 2025}
\DOI{https://doi.org/10.22111/ijfs.2025.50906.8994}

{\let\newpage\relax\vspace*{.1mm}
\noindent\rule{\textwidth}{.3mm}\maketitle}

\thispagestyle{fancylogo}

\mabstract

\section{Introduction}\label{sec_intro}
The theory of fuzzy relational equations (FRE) refers to a generalized version of Boolean relation equations. FRE was originally developed by Sanchez as a formalized model for non-precise concepts \cite{ref_39}. Many applications such as medical diagnosis and treatment are based on fuzzy set theory due to its several properties which make it suitable for formulating uncertain information. Hence, FRE theory was originally applied to medical diagnosis \cite{ref_39}. Considering the sets and various operations, Pedrycz identified two different ways to extend FRE \cite{ref_35}.
Throughout the years, FRE has been utilized in a variety of applications in engineering.
Upon applying the inference rules and identifying the relevant outcomes, the problem of deriving precedence is eased and mathematically reduced to solving an FRE \cite{ref_33}. The fact remains, however, that a large number of issues can potentially be viewed as FRE issues \cite{ref_34}. Due to the fact that FRE holds a number of important applications in various practical fields, many authors have focused their attention on FRE theoretical research, such as its resolution approach and specific optimization problems.

In dealing with FRE problems, it is imperative to identify the solvability of the problem and determine the set of solutions. Based on the study reported in \cite{ref_06}, the authors concluded that the solution set for FRE with a continuous max-t-norm composition (when it is not empty), could be a non-convex set, which is entirely determined with only one maximum solution and a finite number of minimal solutions.
It faces two significant bottlenecks when it comes to tackling FRE-related problems. The first is the non-convexity property of solution sets. This property significantly increases the complexity of problems, especially in the case of optimization problems involving FREs. The second bottleneck lies in efficiently identifying minimal solutions.
As described in Chen and Wang's study \cite{ref_02}, they proposed an algorithm for extracting all minimal solutions and concluded that there might not be a polynomial-time algorithm for finding all minimal solutions of FRE (with max-min composition).
Furthermore, Markovskii demonstrated a close connection between solving max-product FREs and the covering problem, a well-known NP-hard problem \cite{ref_32}.
Moreover, this result extends beyond the minimum and product t-norm operators. Indeed, it generalizes to a wider class of t-norm functions \cite{ref_03,ref_28, ref_29}.
In the past few decades, many researchers have attempted to investigate the solvability of FRE through different max-t compositions \cite{ref_17,ref_19,ref_41,ref_43,ref_44}.
Building on prior research, several scholars have made significant contributions to the field of fuzzy relational inequalities (FRI). This effort has introduced new concepts, improved time complexity, and refined existing theoretical structures and applications \cite{ref_14,ref_18,ref_19,ref_20,ref_25,ref_51,ref_56}.
In a noteworthy contribution to the field, Li and Yang \cite{ref_25} examined FRI using addition-min composition and developed an algorithm to search for minimal solutions. Their work centered on applying FRI to the data transmission mechanism in a BitTorrent-like peer-to-peer file-sharing system.
As noted in \cite{ref_14}, a mixed fuzzy system was studied by combining two fuzzy relational inequalities $A\varphi x\leq b^1$ and $D\varphi x\geq b^2$, where $\varphi$ refers to an operator which has (closed) convex solutions.

Recent studies have demonstrated that the optimization problem involving FRE and FRI is one of the most fascinating and engaging research topics of its kind \cite{ref_01,ref_10,ref_13,ref_14,ref_17,ref_18,ref_19,ref_20,ref_21,ref_22,ref_24,ref_30,ref_42,ref_45,ref_51,ref_63}.
Due to the effectiveness of integer linear programming (ILP) for this class of problem, as a commonly used method, researchers have designed various methods based on translating the main problem into an ILP problem. This mapping enables well-established techniques to efficiently solve the problem.
Alternatively, other algorithms emphasize the resolution of the feasible region and outline some necessary and sufficient conditions for simplification and optimization processes.
The majority of methods in this category leverage analytical results, which are mostly presented by Sanchez \cite{ref_40} and Pedrycz \cite{ref_36}.
As an example, Fang and Li mapped a linear optimization problem involving FRE with max-min composition into an integer programming problem and solved it by a branch-and-bound method using a jump-tracking technique \cite{ref_11}.
In order to improve the methods that were used by Fang and Li, Wu et al. focused on decreasing the search domain, and based on a necessary condition, they derived three rules for the simplification process \cite{ref_46}.
In Chang and Shieh's work \cite{ref_01}, by analyzing the linear optimization problem subject to fuzzy max–min relation equations, they proposed some novel theoretical results that yielded an enhanced upper bound on the optimal objective value, determining some rules for simplifying the problem and suggesting a rule for reducing the solution tree.
As an application of optimizing the linear objective with max-min composition, the author in \cite{ref_23}, implemented a three-tier framework for the streaming media provider to achieve the minimum cost while meeting the requirements.
In other approaches, such as Loetamonphong and Fang, they separated negative from non-negative coefficients in the objective function, dividing the objective function into two subproblems, and by combining the optimal solutions from the two subproblems, an optimal solution was found \cite{ref_31}.
In addition, many authors have studied modifications of the linear optimization problem subject to FRE that were defined with other compositions rather than max-min or max-product, like max-average composition \cite{ref_45} or max-t-norm composition \cite{ref_17,ref_19,ref_21,ref_24,ref_42}.
Accordingly, in order to address the linear optimization problem subjected to a system of sup-t equations, Li and Fang reduced the problem to a 0-1 integer optimization problem \cite{ref_24}.
Authors in \cite{ref_21} outline a very method that deals with linear optimization problems involving fuzzy relation equations with constraints defined with the max-Archimedean t-norm.
It is worth mentioning that a study reported in \cite{ref_42} examined the problem in equal detail with the continuous Archimedean t-norm that employed the covering problem instead of the branch-and-bound methods. 

Over the last few years, researchers in the field of FRE have often investigated the generalization form of linear programming applied to the system of fuzzy relations. They have explored it further by considering the composite operations used in FRE, fuzzy relations used in constraint definition, their impact on the objective function of the problems, and other approaches \cite{ref_05,ref_07,ref_12,ref_19,ref_27,ref_30,ref_48,ref_52}.
In this context, Wu et al. outline an efficient method for the optimization of linear fractional programming problems involving FRE with a max-Archimedean t-norm composition \cite{ref_48}.
Accordingly, Dempe and Ruziyeva aimed at extending the fuzzy linear optimization problem by taking fuzzy coefficients into account \cite{ref_05}.
A further study was conducted by Dubey et al. \cite{ref_07} regarding linear programming problems that utilize interval uncertainty modeled in terms of intuitionist fuzzy sets.
As reported by Yang \cite{ref_52}, the optimal solution to the minimization of a linear objective function examined under the given constraints is referred to as $\sum_{j = 1}^{n} min\left\{ a_{ij}, x_j\right\} \geq b_i $ for $i=1,\dots,m$.
As well, in \cite{ref_51}, the writers discussed the latticized linear programming problem subject to FRI that the constraint is defined as max-product with application to the optimization management model for wireless communication emission base stations.  
The objective function for the latticized linear programming problem was defined as $z(x)=\min_{j=1}^n\{ x_j \}$ within the feasible region $X(A,b)=\left\{ x\in [0,1]^n:A\circ x\geq b \right\}$ that "$\circ$" refers to the fuzzy max-product composition, for which an algorithm was proposed by considering the resolution of the feasible region.

A bipolar fuzzy relation equation is an intriguing generalization of the fuzzy relation equation. Information processing and preference in humans can be viewed as the main example of bipolarity \cite{ref_08}.
It should be noted that Dubois and Prade have reported a study about the relation between asymmetric bipolar and possibility theory \cite{ref_09}. As discussed, possibility theory provides a suitable framework for dealing with bipolar representations, as it allows differentiation between negative and positive information in preference modeling \cite{ref_08,ref_09}.
An examination of a linear optimization problem under bipolar FRE was also the main interest of various authors where FRE was defined with max-min (with application to evaluate public awareness of the products of a supplier) \cite{ref_12,ref_26}, max-product \cite{ref_04,ref_58} and max-Lukasiewicz composition \cite{ref_27,ref_30,ref_50,ref_53}.
The first attempt to address the concept of bipolar FRE was published in \cite{ref_12} which applied max-min composition where constraints are defined as $max_{j=1}^{n}\left\{max \left\{min \left\{a_{i j}^{+}, x_{j}\right\}, min \left\{a_{i j}^{-}, 1-x_{j}\right\}\right\}\right\}$ for $i=1, \dots, m$, where $a_{i j}^{+}, a_{i j}^{-}, x_{j} \in[0,1]$.
Likewise, the writers of \cite{ref_27} outlined a linear optimization problem based upon a bipolar FRE system defined as $X\left(A^{+}, A^{-}, b\right)=\left\{x \in[0,1]^{m}: x \circ A^{+} \vee \tilde{x} \circ A^{-}=b\right\}$ where $\tilde{x}_{i}=1-x_{i}$ for each element of $\tilde{x}=\left(\tilde{x}_{i}\right)_{1 \times m}$ and the notations "$\vee$" and "$\circ$" respectively refer to max operation and the max-Lukasiewicz composition.
The authors resolved the problem by translating the original problem into a 0-1 integer linear problem.
Alternatively, the discussed problem was addressed by an analytical method in \cite{ref_30} that is motivated by the resolution and some structural properties of the feasible region (by applying a necessary condition for the determination of an optimal solution and a simplification process to reduce the problem).
There is, however, no mention of the resolution method used in \cite{ref_50} for determining the complete solution set of the bipolar max-Lukasiewicz FRE.
As a result of Yang's research on bipolar max-Lukasiewicz FRE \cite{ref_50}, it has been demonstrated that the complete solution set of the system is completely specified by a finite number of conservative bipolar paths.
In addition, the authors in \cite{ref_60} studied a more general class of bipolar FREs defined by max strict t-norms.
They categorized the constraints of the problem into four groups, each requiring a specific approach to determine the corresponding feasible solution set.
Moreover, unlike the bipolar path approach \cite{ref_50}, the authors demonstrated that modified ordinary FRE paths could be used to solve the problem, which is achieved by combining the FRE paths specified in the four groups of the constraints \cite{ref_60}.

This paper studies a wide class of bipolar FRE linear optimization problems as the following mathematical model in which $\varphi$ is an arbitrary continuous Archimedean t-norm:

\begin{equation}\label{eq_1}
	\begin{array}{ll}
		\min & c^{T} x                                                   \\
		& A^{+} \varphi x \vee A^{-} \varphi(\mathbf{1}-x) = b \\
		& x \in[0,1]^{n}
	\end{array}
\end{equation}

where $A^{+}=\left(a_{i j}^{+}\right)_{m \times n}$ and $A^{-}=\left(a_{i j}^{-}\right)_{m \times n}$ are fuzzy matrices and $b=\left(b_{i}\right)_{m \times 1}$ is a fuzzy vector such that $0 \leq a_{i j}^{+} \leq 1,0 \leq a_{i j}^{-} \leq 1$ and $0 \leq b_{i} \leq 1$ for each $i \in \mathscr{I}=\{1,2, \ldots, m\}$ and each $j \in \mathscr{J}=\{1,2, \ldots, n\}$, respectively. As well, the constant vector $c=\left(c_{j}\right)_{n \times 1}$, the sum vector \textbf{1} (each component of the vector is equal to one) and the unknown vector $x=\left(x_{j}\right)_{n \times 1}$ are in $\mathbb{R} ^{n}$. Furthermore, $\varphi$ refers to an arbitrary continuous Archimedean t-norm, whereas in Problem (\ref{eq_1}) the constraints are defined as:  

\begin{equation}\label{eq-2}
	\max _{j=1}^{n}\left\{\max \left\{\varphi\left(a_{i j}^{+}, x_{j}\right), \varphi\left(a_{i j}^{-}, 1-x_{j}\right)\right\}\right\}=b_{i} \quad, \forall i \in \mathscr{I}
\end{equation}

It is worth mentioning that since the objective function can be transformed into an increasing function using the change of variables \cite{ref_12}, it can be assumed without loss of generality that $c_{j} \geq 0$, $\forall i \in \mathscr{I}$. In the current method, unlike the method presented in \cite{ref_60}, to find a global optimum, it is not necessary to obtain all local optimal solutions.
Moreover, by using the current method, the constraints do not need to be classified (and therefore, it is unnecessary to apply different approaches to check the feasibility of constraints in different groups and the optimality of solutions in their feasible regions).
It is proved that the feasible region of Problem (\ref{eq_1}) can be determined as a finite number of compact sets, which are obtained by functions called admissible functions in this paper. After deriving some basic properties of Archimedean bipolar fuzzy equations, seven simplification techniques are presented which are applied at the beginning of the algorithm to reduce the problem size; four techniques are based on eliminating redundant constraints, and three additional techniques eliminate a set of feasible solutions that are not optimal.
Subsequently, after applying the simplification techniques, the proposed algorithm takes advantage of two accelerating strategies to find an optimal solution. The first strategy removes some candidate solutions that are not optimal (by reducing the domain search to a subset of admissible functions called modified functions in this paper), and the second strategy (a branch and bound method), removes some branches that do not lead to an optimal solution.

The rest of the paper is organized as follows. Section \ref{sec-2} presents preliminary definitions, concepts, and properties of Archimedean bipolar FREs. In Section \ref{sec-3}, the feasible solution set of the original problem is characterized, and two necessary conditions are derived to identify the feasibility of the problem. Section \ref{sec-4} describes a necessary optimality condition that reduces the search domain by focusing on a subset of the feasible region rather than the entire set of feasible solutions. In Section \ref{sec-5}, seven simplification techniques are introduced to accelerate the resolution process by reducing the problem size. By taking advantage of the concept of modified functions, a branch-and-bound algorithm is presented to solve Problem (\ref{eq_1}), and, finally in Section \ref{sec-6} an illustration of this algorithm is provided using an example.
\section{Preliminary definitions and properties}\label{sec-2}
This section initially describes the feasible solution set of the equation $\varphi(a, x)=b$ where $a$ and $b$ are two fixed scalars in [0,1], $x \in[0,1]$ and $\varphi$ is an arbitrary continuous Archimedean t-norm. Then, the feasible solution set of the equation $\max \left\{\varphi\left(a_{i j}^{+}, x_{j}\right), \varphi\left(a_{i j}^{-}, 1-x_{j}\right)\right\}=b_{i}$ is completely characterized for each $i \in \mathscr{I}$ and each $j \in \mathscr{J}$. Subsequently, according to the results obtained in this section, the feasible region of Problem (\ref{eq_1}) is determined in the next section. For the sake of simplicity, let $S_{i}$ denote the feasible solution set of the $i$'th equation; that is, $S_{i}=\left\{x \in[0,1]^{n}: \max _{j=1}^{n}\left\{\max \left\{\varphi\left(a_{i j}^{+}, x_{j}\right), \varphi\left(a_{i j}^{-}, 1-x_{j}\right)\right\}\right\}=b_{i}\right\}$. Also, let $S\left(A^{+}, A^{-}, b\right)$ denote the feasible solution set for Problem (\ref{eq_1}). So, it is clear that $S\left(A^{+}, A^{-}, b\right)=\bigcap_{i \in \mathscr{I}} S_{i}$.
\begin{definition}\label{def-1}
	Suppose that $a$ and $b$ are two fixed scalers in [0,1]. We define $S=\{x \in[0,1]: \varphi(a, x)=b\}$ and $I=\{x \in[0,1]: \varphi(a, x) \leq b\}$. Also, if $S \neq \varnothing$, let $l=\inf S$ and $u=\sup S$.
\end{definition}

It should be noted that if $S \neq \varnothing$ in Definition \ref{def-1}, then both $l=\inf S$ and $u=\sup S$ exist. This consequence is actually resulted from the fact that $S \subseteq[0,1]$ and the least-upper-bound property of the real numbers set.
\begin{lemma}\label{lm-1}
	Suppose that $\varphi$ is a continuous Archimedean t-norm,
	\textbf{(a)} If $x_{0} \in S$, then $x \in I$ for each $x \in\left[0, x_{0}\right]$. \textbf{(b)} If $a<b$, then $I=[0,1]$ and $S=\varnothing$.
\end{lemma}

\begin{proof}
	\textbf{(a)} Since $x_{0} \in S$, then from Definition \ref{def-1} we have $\varphi\left(a, x_{0}\right)=b$. On the other hand, from the monotonicity property of t-norms, $\varphi(a, x) \leq \varphi\left(a, x_{0}\right)$, $\forall x \in\left[0, x_{0}\right]$. So, $\varphi(a, x) \leq b$, $\forall x \in\left[0, x_{0}\right]$, which together with Definition \ref{def-1} imply $x \in I$, $\forall x \in\left[0, x_{0}\right]$. \textbf{(b)} The result follows from the identity and monotonicity laws of t-norms and the inequalities $\varphi(a, x) \leq \varphi(a, 1)=a<b$, $\forall x \in[0,1]$.
\end{proof}

\begin{lemma}{\cite{ref_61}}\label{lm-2}
	A function $\varphi:[0,1]^{2} \rightarrow[0,1]$ is a continuous Archimedean t-norm if and only if there exist a strictly decreasing and continuous function $f_{\varphi}:[0,1] \rightarrow[0,+\infty)$ with $f_{\varphi}(1)=0$ such that
	\begin{equation}\label{eq-3}
		\varphi(x, y)=f_{\varphi}^{(-1)}\left(f_{\varphi}(x)+f_{\varphi}(y)\right)
	\end{equation}
	where $f_{\varphi}^{(-1)}$ is the pseudoinverse of $f_{\varphi}$ defined by
	\begin{equation}\label{eq-4}
		f_{\varphi}^{(-1)}(x)= \begin{cases}f_{\varphi}^{-1}(x) & \text { if } x \leq f_{\varphi}(0) \\ 0 & \text { otherwise }\end{cases}
	\end{equation}
	and $f^{-1}_{\varphi}(x)$ is the inverse of $f_{\varphi}$. Moreover, representation (\ref{eq-3}) is unique up to a positive multiplicative constant.
	_\end{lemma}

We say that $\varphi$ is generated by $f_{\varphi}$ if $\varphi$ has representation (\ref{eq-3}). In this case, $f_{\varphi}$ is said to be an additive generator of $\varphi$. For example, an additive generator of Product tnorm is given by $f_{P}(x)=-\operatorname{Ln}(x)$. Then, $f_{P}^{(-1)}(x)=e^{-x}$. Also, Lukasiewicz t-norm is generated by $f_{L}(x)=1-x$ with $f_{L}^{(-1)}(x)=\max \{1-x, 0\}$. Table \ref{tbl_a1} (see Appendix A) lists some frequently used continuous Archimedean t-norms (including strict and nilpotent $\mathrm{t}$-norms). For each $\mathrm{t}$-norm in Table \ref{tbl_a1}, its additive generator $f_{\varphi}$ and pseudoinverse $f_{\varphi}^{(-1)}$ discussed in Lemma \ref{lm-2} are presented in Table \ref{tbl_a2}.
\begin{remark}\label{rmk-1}
	Suppose that $a \in[0,1]$ and $0<b \leq 1$. Based on Lemma \ref{lm-2}, the equation $\varphi(a, x)=b$ can be rewritten as follows:
	\begin{equation}\label{eq-5}
		\varphi(a, x)=f_{\varphi}^{(-1)}\left(f_{\varphi}(a)+f_{\varphi}(x)\right)=f_{\varphi}^{-1}\left(f_{\varphi}(a)+f_{\varphi}(x)\right)=b
	\end{equation}
	Also, if $b=0$ we have
	\begin{equation}\label{eq-6}
		\varphi(a, x)=f_{\varphi}^{(-1)}\left(f_{\varphi}(a)+f_{\varphi}(x)\right)=b=0
	\end{equation}
\end{remark}

\begin{theorem}\label{thm-1}
	Suppose that $\varphi$ is a continuous Archimedean t-norm, $f_{\varphi}$ is an additive generator of $\varphi$ and $u=f_{\varphi}^{-1}(f_{\varphi}(b)-$ $f_{\varphi}(a))$, where $a \geq b>0$. \textbf{(a)} $S=\{u\}$ and $I=[0, u]$. \textbf{(b)} If $b<a<1$, then $b<u<1$. \textbf{(c)} If $b<a=1$, then $u=b$. \textbf{(d)} If $a=b$, then $u=1$.
\end{theorem}

\begin{proof}
	\textbf{(a)} Clearly From Remark \ref{rmk-1}, $\varphi(a,u)=f_{\varphi}^{-1}(f_{\varphi}(a)+f_{\varphi}(u))=f_{\varphi}^{-1}(f_{\varphi}(a)+(f_{\varphi}(b)-f_{\varphi}(a)))$. Also relation (\ref{eq-5}) imply $\varphi(a, x)=b$ if and only if $f_{\varphi}^{-1}\left(f_{\varphi}(a)+f_{\varphi}(x)\right)=b$. Now, since $f_{\varphi}$ and $f_{\varphi}^{-1}$ are strictly decreasing continuous functions, it is concluded that $f_{\varphi}^{-1}\left(f_{\varphi}(a)+f_{\varphi}(x)\right)=b$ if and only if $x=u=f_{\varphi}^{-1}\left(f_{\varphi}(b)-f_{\varphi}(a)\right)$. Hence, $S=\{u\}$, which together with Lemma \ref{lm-1} (a) imply $I=[0, u]$. \textbf{(b)} From Part (a), we have $S=\{u\}$. Since $\varphi(a, 1)=a>b$, then $1 \notin S$, which implies $u<1$. On the other hand, for each $x \in[0, b)$, $\varphi(a, x) \leq \varphi(1, x)=x<b$, that means $x \notin S$, $\forall x \in[0, b)$. Moreover, $b \notin S$; because otherwise, if $\varphi(a, b)=b$, then the associativity property of t-norms implies $\varphi(a, \varphi(a, b))=\varphi(\varphi(a, a), b)=b$, and by repeating the argument we have $\varphi\left(\varphi^{(n)}(a), b\right)=b$ ($\forall n \in \mathbb{R} ^{n}$), where $\varphi^{(n)}(a)$ denotes $\varphi(\underbrace{a, \ldots, a}_{n \text {-times }})$. Therefore,
	\begin{equation}\label{eq-7}
		\lim _{n \rightarrow \infty} \varphi\left(\varphi^{(n)}(a), b\right)=b
	\end{equation}
	But, since $\varphi$ is continuous and Archimedean, it is concluded that
	\begin{equation}\label{eq-8}
		\lim _{n \rightarrow \infty} \varphi\left(\varphi^{(n)}(a), b\right)=\varphi\left(\lim _{n \rightarrow \infty} \varphi^{(n)}(a), b\right)=\varphi(0, b)=0
	\end{equation}
	Hence, from (\ref{eq-7}) and (\ref{eq-8}), $b=0$ that is a contradiction. Consequently, for each $x \in[0, b]$, we have $x \notin S$, that results in $b<u$. \textbf{(c)} In this case, the equation $\varphi(a, x)=b$ is converted to $\varphi(1, x)=x=b$, that is, $S=\{b\}$. So, the result follows from $S=\{b\}$ and Part (a). \textbf{(d)} At first, we note that if $a=b$, then $\varphi(a, 1)=a=b$; that is, $1 \in S$, which together with Part (a) imply $u=1$.
\end{proof}

\begin{theorem}\label{thm-2}
	Suppose that $\varphi$ is a continuous Archimedean t-norm, $f_{\varphi}$ is an additive generator of $\varphi$ and $u=f_{\varphi}^{-1}(f_{\varphi}(b)-$ $f_{\varphi}(a))$, where $a \geq b=0$. \textbf{(a)} If $a=b=0$, then $S=I=[0, u]$ and $u=1$. \textbf{(b)} If $a>b=0$, then $S=I=[0, u]$ and $0 \leq u<1$.
\end{theorem}

\begin{proof}
	\textbf{(a)} If $a=b=0$, then the equation $\varphi(a, x)=b$ is converted to $\varphi(0, x)=0$, which is satisfied for each $x \in[0,1]$. Hence, $S=[0,1]$ and therefore $I=[0,1]$ from Lemma \ref{lm-1} (a). Also, $u=f_{\varphi}^{-1}\left(f_{\varphi}(b)-f_{\varphi}(a)\right)=f_{\varphi}^{-1}\left(f_{\varphi}(0)-f_{\varphi}(0)\right)=f_{\varphi}^{-1}(0)=1$. \textbf{(b)} In this case, from relation (\ref{eq-6}), the equation $\varphi(a, x)=0$ is converted to $f_{\varphi}^{(-1)}\left(f_{\varphi}(a)+f_{\varphi}(x)\right)=0$, which together with (\ref{eq-4}) imply
	\begin{equation}\label{eq-9}
		f_{\varphi}(a)+f_{\varphi}(x)>f_{\varphi}(0)
	\end{equation}
	So, $f_{\varphi}(x)>f_{\varphi}(0)-f_{\varphi}(a)$ that requires $x<f_{\varphi}^{-1}\left(f_{\varphi}(0)-f_{\varphi}(a)\right)=u$. Hence, we have
	\begin{equation}\label{eq-10}
		\varphi(a, x)=0, \forall x \in[0, u)
	\end{equation}
	Based on (\ref{eq-10}), it follows that $\lim _{x \rightarrow u^{-}} \varphi(a, x)=0$. Furthermore, since $\varphi$ is continuous,
	$$\lim _{x \rightarrow u^{-}} \varphi(a, x)=\varphi\left(a, \lim _{x \rightarrow u^{-}} x\right)=\varphi(a, u).$$
	So, $\varphi(a, u)=0$ which together with (\ref{eq-10}) imply $\varphi(a, x)=0$, $\forall x \in[0, u]$; that is, $[0, u] \subseteq S$. In order to prove $S=[0, u]$, we shall show that $x \notin S$ for each $x>u$. By contradiction, suppose that $\varphi(a, x)=0$ (i.e., $x \in S$ ) and $x>u=f_{\varphi}^{-1}\left(f_{\varphi}(0)-f_{\varphi}(a)\right)$. Since $\varphi(a, x)=0$, then $x$ satisfies (\ref{eq-9}). On the other hand, since $x>f_{\varphi}^{-1}\left(f_{\varphi}(0)-f_{\varphi}(a)\right)$, it is concluded that $f_{\varphi}(x)<f_{\varphi}(0)-f_{\varphi}(a)$, that contradicts (\ref{eq-4}). Consequently, $S=[0, u]$ and therefore $I=[0, u]$ from Lemma \ref{lm-1} (a). Moreover, since $\varphi(a, 1)=a>b$, it follows that $1 \notin S=[0, u]$, i.e., $u<1$.
\end{proof}

From Lemma \ref{lm-1} (b), Theorem \ref{thm-1}, and Theorem \ref{thm-2}, we can obtain a complete characterization of $S$ when $\varphi$ is a continuous Archimedean t-norm (i.e., $\varphi$ is either continuous strict or continuous nilpotent). These results are summarized in the following corollary.

\begin{corollary}\label{corl-1}
	Suppose that $\varphi$ is a continuous Archimedean t-norm, $f_{\varphi}$ is an additive generator of $\varphi$ and $u=f_{\varphi}^{-1}(f_{\varphi}(b)-$ $f_{\varphi}(a))$. Then,
	$$
	S=\begin{cases}
		\{u\}       & , a \geq b>0 \\
		{[0, u]}    & , a \geq b=0 \\
		\varnothing & , a<b
	\end{cases}
	\quad,\quad
	I=\begin{cases}
		{[0, u]} & , a \geq b \geq 0 \\
		{[0,1]}  & , a<b
	\end{cases}
	$$
	where $b<u<1$ if $0<b<a<1$; $u=b$ if $0<b<a=1$; $u=1$ if $a=b>0$; $u<1$ if $a>b=0$; and $u=1$, if $a=b=0$.
\end{corollary}

It should be noted that if $a>b=0$, we generally have $u \geq 0$. For example, if $\varphi$ is Lukasiewicz t-norm and $a>b=0$, then $u=1-a$ which is positive if $a<1$. However, if $\varphi$ is strict and $a>b=0$, then for each $x>0$, we have $\varphi(a, x)>\varphi(a, 0)=0$, i.e., $S=\{0\}$. For each t-norm in Table \ref{tbl_a1} (see Appendix A), value $u$ discussed in Theorem \ref{thm-1} and Theorem \ref{thm-2} is presented in Table \ref{tbl_a3}.

\begin{definition}\label{def-2}
	For each $i \in \mathscr{I}$ and each $j \in \mathscr{J}$, we define $$S_{i j}^{+}=\left\{x_{j} \in[0,1]: \varphi\left(a_{i j}^{+}, x_{j}\right)=b_{i}\right\} ~ \text{and} ~ S_{i j}^{-}=\left\{x_{j} \in[0,1]: \varphi\left(a_{i j}^{-}, 1-x_{j}\right)=b_{i}\right\};$$ that is, $S_{i j}^{+}$ and $S_{i j}^{-}$ denote the feasible solution sets of the equations $\varphi\left(a_{i j}^{+}, x\right)=b_{i}$ and $\varphi\left(a_{i j}^{-}, 1-x\right)=b_{i}$, respectively. Furthermore, let $I_{i j}^{+}=\left\{x_{j} \in[0,1]: \varphi\left(a_{i j}^{+}, x_{j}\right) \leq b_{i}\right\}$ and $I_{i j}^{-}=\left\{x_{j} \in[0,1]: \varphi\left(a_{i j}^{-}, 1-x_{j}\right) \leq b_{i}\right\}$.
\end{definition}

According to Definition \ref{def-2}, it is easy to verify that $S_{i j}^{+} \subseteq I_{i j}^{+}$ and $S_{i j}^{-} \subseteq I_{i j}^{-}$, $\forall i \in \mathscr{I}$ and $\forall j \in \mathscr{J}$.

\begin{corollary}\label{corl-2}
	Suppose that $\varphi$ is a continuous Archimedean t-norm and $f_{\varphi}$ is an additive generator of $\varphi$. Also, for each $i \in \mathscr{I}$ and $j \in \mathscr{J}$, let $u_{i j}^{+}=f_{\varphi}^{-1}\left(f_{\varphi}\left(b_{i}\right)-f_{\varphi}\left(a_{i j}^{+}\right)\right)$ and $u_{i j}^{-}=f_{\varphi}^{-1}\left(f_{\varphi}\left(b_{i}\right)-f_{\varphi}\left(a_{i j}^{-}\right)\right)$.
	\begin{itemize}
		\item[\textbf{(a)}]
		If $a_{i j}^{+}<b_{i}\left(a_{i j}^{-}<b_{i}\right)$, then $I_{i j}^{+}=[0,1]$ and $S_{i j}^{+}=\varnothing\left(I_{i j}^{-}=[0,1]\right.$ and $\left.S_{i j}^{-}=\varnothing\right)$.
		\item[\textbf{(b)}]
		If $a_{i j}^{+} \geq b_{i}>0$, then $S_{i j}^{+}=\left\{u_{i j}^{+}\right\}$ and $I_{i j}^{+}=\left[0, u_{i j}^{+}\right]$. Moreover, $b_{i}<u_{i j}^{+}<1$ if $b_{i}<a_{i j}^{+}<1$; $u_{i j}^{+}=b_{i}$ if $b_{i}<a_{i j}^{+}=1$; and $u_{i j}^{+}=1$ if $a_{i j}^{+}=b_{i}$.
		\item[\textbf{(c)}]
		If $a_{i j}^{+} \geq b_{i}=0$, then $S_{i j}^{+}=I_{i j}^{+}=\left[0, u_{i j}^{+}\right]$. Moreover, $u_{i j}^{+}=1$ if $a_{i j}^{+}=b_{i}=0$; and $0 \leq u_{i j}^{+}<1$ if $a_{i j}^{+}>b_{i}=0$.
		\item[\textbf{(d)}]
		If $a_{i j}^{-} \geq b_{i}>0$, then $S_{i j}^{-}=\left\{1-u_{i j}^{-}\right\}$ and $I_{i j}^{-}=\left[1-u_{i j}^{-}, 1\right]$. Moreover, $b_{i}<u_{i j}^{-}<1$ if $b_{i}<a_{i j}^{-}<1$; $u_{i j}^{-}=b_{i}$ if $b_{i}<a_{i j}^{-}=1$; and $u_{i j}^{-}=1$ if $a_{i j}^{-}=b_{i}$.
		\item[\textbf{(e)}]
		If $a_{i j}^{-} \geq b_{i}=0$, then $S_{i j}^{-}=I_{i j}^{-}=\left[1-u_{i j}^{-}, 1\right]$. Moreover, $u_{i j}^{-}=1$ if $a_{i j}^{-}=b_{i}=0$; and $0 \leq u_{i j}^{-}<1$ if $a_{i j}^{-}>b_{i}=0$.
	\end{itemize}
\end{corollary}

\begin{proof}
	The results directly follow from Lemma \ref{lm-1} (b), Theorem \ref{thm-1}, Theorem \ref{thm-2}, Definition \ref{def-2} and replacing $a_{i j}^{+}$ and $a_{i j}^{-}$ by $a$; $b_i$ by $b$; $S_{i j}^{+}$ and $S_{i j}^{-}$ by $S$; and $I_{i j}^{+}$ and $I_{i j}^{-}$ by $I$.
\end{proof}

\begin{definition}\label{def-3}
	For each $i \in \mathscr{I}$ and each $j \in \mathscr{J}$, we define $S_{i j}=\left\{x_{j} \in[0,1]: \max \left\{\varphi\left(a_{i j}^{+}, x_{j}\right), \varphi\left(a_{i j}^{-}, 1-x_{j}\right)\right\}=b_{i}\right\}$ and $I_{i j}=\left\{x_{j} \in[0,1]: \max \left\{\varphi\left(a_{i j}^{+}, x_{j}\right), \varphi\left(a_{i j}^{-}, 1-x_{j}\right)\right\} \leq b_{i}\right\}$.
\end{definition}

\begin{lemma}\label{lm-3}
	\textbf{(a)} $I_{i j}=I_{i j}^{+} \cap I_{i j}^{-}$. \textbf{(b)} $S_{i j}=I_{i j} \cap\left(S_{i j}^{+} \cup S_{i j}^{-}\right)$.
\end{lemma}

\begin{proof}
	\textbf{(a)} The proof is easily resulted from Definitions 2 and 3. \textbf{(b)} According to Definition \ref{def-3}, $x \in S_{i j}$ if and only if $\varphi\left(a_{i j}^{+}, x_{j}\right) \leq b_{i}$, $\varphi\left(a_{i j}^{-}, 1-x_{j}\right) \leq b_{i}$ (i.e., $x \in I_{i j}^{+} \cap I_{i j}^{-}=I_{i j}$) and at least one of the two equalities $\varphi\left(a_{i j}^{+}, x_{j}\right)=b_{i}$ (i.e., $\left.x \in S_{i j}^{+}\right)$ and $\left.\varphi\left(a_{i j}^{-}, 1-x_{j}\right)\right\}=b_{i}$ (i.e., $x \in S_{i j}^{-}$) holds.
\end{proof}

Based on Corollary \ref{corl-2} and Lemma \ref{lm-3}, we can obtain the following results, which characterize sets $S_{i j}$ and $I_{i j}$ for all cases.

\begin{corollary}\label{corl-3}
	Suppose that $\varphi$ is a continuous Archimedean t-norm and $f_{\varphi}$ is an additive generator of $\varphi$. Also, for each $i \in \mathscr{I}$ and $j \in \mathscr{J}$, let $u_{i j}^{+}=f_{\varphi}^{-1}\left(f_{\varphi}\left(b_{i}\right)-f_{\varphi}\left(a_{i j}^{+}\right)\right)$ and $u_{i j}^{-}=f_{\varphi}^{-1}\left(f_{\varphi}\left(b_{i}\right)-f_{\varphi}\left(a_{i j}^{-}\right)\right)$.
	\begin{itemize}
		\item[\textbf{(a)}]
		If $a_{i j}^{+}<b_{i}$ and $a_{i j}^{-}<b_{i}$, then $I_{i j}=[0,1]$ and $S_{i j}=\varnothing$.
		\item[\textbf{(b)}]
		If $a_{i j}^{+} \geq b_{i}>0$ and $a_{i j}^{-}<b_{i}$, then $S_{i j}=\left\{u_{i j}^{+}\right\}$ and $I_{i j}=\left[0, u_{i j}^{+}\right]$.
		\item[\textbf{(c)}]
		If $a_{i j}^{-} \geq b_{i}>0$ and $a_{i j}^{+}<b_{i}$, then $S_{i j}=\left\{1-u_{i j}^{-}\right\}$ and $I_{i j}=\left[1-u_{i j}^{-}, 1\right]$.
		\item[\textbf{(d)}]
		If $a_{i j}^{+} \geq b_{i}=0$ and $a_{i j}^{-} \geq b_{i}=0$, then $S_{i j}=I_{i j}=\left[1-u_{i j}^{-}, u_{i j}^{+}\right]$.
		\item[\textbf{(e)}]
		If $a_{i j}^{+} \geq b_{i}>0$ and $a_{i j}^{-} \geq b_{i}>0$, then $S_{i j}=\left\{1-u_{i j}^{-}, u_{i j}^{+}\right\}$ and $I_{i j}=\left[1-u_{i j}^{-}, u_{i j}^{+}\right]$.
	\end{itemize}
\end{corollary}

The results of the above-mentioned corollaries can be summarized as follows:

\begin{remark}\label{rmk-2}
	Suppose that $\varphi$ is a continuous Archimedean t-norm and $f_{\varphi}$ is an additive generator of $\varphi$. Also, for each $i \in \mathscr{I}$ and $j \in \mathscr{J}$, let $u_{i j}^{+}=f_{\varphi}^{-1}\left(f_{\varphi}\left(b_{i}\right)-f_{\varphi}\left(a_{i j}^{+}\right)\right)$ and $u_{i j}^{-}=f_{\varphi}^{-1}\left(f_{\varphi}\left(b_{i}\right)-f_{\varphi}\left(a_{i j}^{-}\right)\right)$. Then,
	$$
	S_{i j}=\left\{\begin{array}{ll}
		\varnothing                               & , a_{i j}^{+}, a_{i j}^{-}<b_{i}                           \\
		\left\{u_{i j}^{+}\right\}                & , a_{i j}^{+} \geq b_{i}>0 \text { and } a_{i j}^{-}<b_{i} \\
		\left\{1-u_{i j}^{-}\right\}              & , a_{i j}^{-} \geq b_{i}>0 \text { and } a_{i j}^{+}<b_{i} \\
		{\left[1-u_{i j}^{-}, u_{i j}^{+}\right]} & , a_{i j}^{+}, a_{i j}^{-} \geq b_{i}=0                    \\
		\left\{1-u_{i j}^{-}, u_{i j}^{+}\right\} & , a_{i j}^{+}, a_{i j}^{-} \geq b_{i}>0
	\end{array} \quad, \quad I_{i j}= \begin{cases}{[0,1]}                                   & , a_{i j}^{+}, a_{i j}^{-}<b_{i}                           \\
		{\left[0, u_{i j}^{+}\right]}             & , a_{i j}^{+} \geq b_{i}>0 \text { and } a_{i j}^{-}<b_{i} \\
		{\left[1-u_{i j}^{-}, 1\right]}           & , a_{i j}^{-} \geq b_{i}>0 \text { and } a_{i j}^{+}<b_{i} \\
		{\left[1-u_{i j}^{-}, u_{i j}^{+}\right]} & , a_{i j}^{+}, a_{i j}^{-} \geq b_{i} \geq 0\end{cases}\right.
	$$
	where $b_{i}<u_{i j}^{+}<1$ if $0<b_{i}<a_{i j}^{+}<1$; $u_{i j}^{+}=b_{i}$ if $0<b_{i}<a_{i j}^{+}=1$; $u_{i j}^{+}=1$ if $a_{i j}^{+}=b_{i}$ ; $0 \leq u_{i j}^{+}<1$ if $a_{i j}^{+}>b_{i}=0$; $b_{i}<u_{i j}^{-}<1$ if $0<b_{i}<a_{i j}^{-}<1$; $u_{i j}^{-}=b_{i}$ if $0<b_{i}<a_{i j}^{-}=1$; $u_{i j}^{-}=1$ if $a_{i j}^{-}=b_{i}$; and $0 \leq u_{i j}^{-}<1$ if $a_{i j}^{-}>b_{i}=0$.
\end{remark}

\begin{remark}\label{rmk-3}
	According to Remark \ref{rmk-2}, each $I_{i j}$ can be written as an interval $I_{i j}=\left[L_{i j}, U_{i j}\right]$ where $L_{i j}$ and $U_{i j}$ are the lower and upper bounds of the interval, respectively. Also, each $S_{i j}$ (if it is non-empty) is obtained as one of the forms of $S_{i j}=\left\{L_{i j}\right\}$, $S_{i j}=\left\{U_{i j}\right\}$, $S_{i j}=\left\{L_{i j}, U_{i j}\right\}$ or $S_{i j}=\left[L_{i j}, U_{i j}\right]$.
\end{remark}

\begin{definition}\label{def-4}
	For each $j \in \mathscr{J}$, we define $I_{j}=\bigcap_{i \in \mathscr{I}} I_{i j}$. Also, let $S_{i j}^{\prime}=S_{i j} \cap I_{j}$, $\forall i \in \mathscr{I}$ and $\forall j \in \mathscr{J}$.
\end{definition}

\begin{corollary}\label{corl-4}
	Let $i \in \mathscr{I}$ and $j \in \mathscr{J}$. If $S_{i j}^{\prime} \neq \varnothing$, then $S_{i j}^{\prime}$ is as one of the forms of $S_{i j}^{\prime}=\left\{L_{j}\right\}$, $S_{i j}^{\prime}=\left\{U_{j}\right\}$, $S_{i j}^{\prime}=\left\{L_{j}, U_{j}\right\}$ or $S_{i j}^{\prime}=\left[L_{j}, U_{j}\right]$, where $L_{j}$ and $U_{j}$ are the lower and upper bound of $I_{j}=\left[L_{j}, U_{j}\right]$, respectively.
\end{corollary}

\begin{proof}
	From Definition \ref{def-4} and Remark \ref{rmk-3}, we have
	$$
	I_{j}=\bigcap_{i \in \mathscr{I}} I_{i j}=\bigcap_{i \in \mathscr{I}}\left[L_{i j}, U_{i j}\right]=\left[\max _{i \in \mathscr{I}}\left\{L_{i j}\right\}, \min _{i \in \mathscr{I}}\left\{U_{i j}\right\}\right]=\left[L_{j}, U_{j}\right]
	$$
	Hence, for each $i \in \mathscr{I}$, if $L_{i j} \in I_{j}$ $\left(U_{i j} \in I_{j}\right)$, then we must have $L_{i j}=L_{j}$ $\left(U_{i j}=U_{j}\right)$. On the other hand, from Remark \ref{rmk-3}, a non-empty $S_{i j}$ is as one of the forms of $S_{i j}=\left\{L_{i j}\right\}$, $S_{i j}=\left\{U_{i j}\right\}$, $S_{i j}=\left\{L_{i j}, U_{i j}\right\}$ or $S_{i j}=\left[L_{i j}, U_{i j}\right]$. So, if $S_{i j}=\left\{L_{i j}\right\}$ and $S_{i j}^{\prime}=S_{i j} \cap I_{j} \neq \varnothing$, then $L_{i j} \in I_{j}$; that is, $L_{i j}=L_{j}$ and therefore $S_{i j}^{\prime}=S_{i j} \cap I_{j}=\left\{L_{j}\right\}$. For the other cases, the results are attained by the similar argument.
\end{proof}

\begin{example}\label{ex-1}
	Consider Problem (\ref{eq_1}) includes Yager operator where $$z=c^{T} x=x_{1}+0.35 x_{2}+0.93 x_{3}+3.28 x_{4}+5.03 x_{5}+2.96 x_{6}+x_{7}+2.75 x_{8}+5.25 x_{9}+6.39 x_{10}$$
	
	$$A^{+}=\left[\begin{array}{llllllllll}
		0.25 & 0.32 & 0.41 & 0.19 & 0.70 & 0.13 & 0.44 & 0.37 & 0.28 & 0.50 \\
		0.80 & 0.73 & 0.64 & 0.79 & 0.80 & 0.22 & 0.80 & 0.56 & 0.10 & 0.28 \\
		0.11 & 0.20 & 0.12 & 0.13 & 0.05 & 0.25 & 0.40 & 0.25 & 0.20 & 0.18 \\
		0.10 & 0.23 & 0.25 & 0.15 & 0.12 & 0.05 & 0.02 & 0.01 & 0.15 & 0.15 \\
		0.45 & 0.35 & 0.70 & 0.50 & 0.41 & 0.27 & 0.39 & 0.48 & 0.17 & 0.39 \\
		0.60 & 0.70 & 0.25 & 0.38 & 0.63 & 0.58 & 0.46 & 0.47 & 0.85 & 0.33 \\
		0.01 & 0.02 & 0.15 & 0.09 & 0.12 & 0.10 & 0.15 & 0.15 & 0.04 & 0.25 \\
		0.75 & 0.64 & 0.32 & 0.29 & 0.39 & 0.61 & 0.57 & 0.34 & 1.00 & 0.46 \\
		0.22 & 0.20 & 0.35 & 0.23 & 0.30 & 0.18 & 0.29 & 0.25 & 0.35 & 0.10 \\
		0.41 & 0.25 & 0.50 & 0.20 & 0.56 & 0.60 & 0.59 & 0.60 & 0.47 & 0.31
	\end{array}\right]$$
	
	$$A^{-}=\left[\begin{array}{llllllllll}
		0.70 & 0.70 & 0.32 & 0.44 & 0.00 & 0.16 & 0.20 & 0.50 & 0.40 & 0.39 \\
		0.70 & 0.65 & 0.14 & 0.12 & 0.80 & 0.76 & 0.00 & 1.00 & 0.15 & 0.79 \\
		0.17 & 0.24 & 0.20 & 0.20 & 0.06 & 0.25 & 0.13 & 0.19 & 0.22 & 0.02 \\
		0.14 & 0.10 & 0.04 & 0.00 & 0.10 & 0.00 & 0.14 & 0.02 & 0.15 & 0.08 \\
		0.70 & 0.04 & 0.27 & 0.36 & 0.60 & 0.40 & 0.48 & 0.50 & 0.50 & 0.50 \\
		0.66 & 0.63 & 0.14 & 0.73 & 0.53 & 0.46 & 0.61 & 0.85 & 0.85 & 0.39 \\
		0.00 & 0.15 & 0.15 & 0.05 & 0.02 & 0.03 & 0.10 & 0.12 & 0.18 & 0.09 \\
		0.63 & 0.03 & 0.55 & 0.77 & 0.79 & 0.49 & 0.21 & 0.32 & 0.80 & 0.71 \\
		0.27 & 0.30 & 0.35 & 0.24 & 0.35 & 0.07 & 0.29 & 0.35 & 0.20 & 0.75 \\
		0.59 & 0.34 & 0.26 & 0.38 & 0.02 & 0.60 & 0.52 & 0.43 & 0.27 & 0.44
	\end{array}\right]$$
	
	$$b^{T}=[0.50,0.80,0.25,0.15,0.50,0.75,0.15,0.80,0.35,0.60]$$
	
	In this example, $\mathscr{I}=\{1,2, \ldots, 10\}$, $\mathscr{J}=\{1,2, \ldots, 10\}$ and $\varphi$ is the Yager t-norm (that is continuous and nilpotent) with $p=2$. So, according to Table \ref{tbl_a1} we have
	$$\varphi(x, y)=T_{Y}^{2}(x, y)=\max \left\{0,1-\left((1-x)^{2}+(1-y)^{2}\right)^{1 / 2}\right\}$$
	Also, $f_{Y}(x)=(1-x)^{2}$ and $f_{Y}^{(-1)}(x)=\max \{1-\sqrt{x}, 0\}$ (see Table \ref{tbl_a2}). For $i=8$ and $j=9$, we have $a_{89}^{+}=1$ and $a_{89}^{-}=b_{8}=0.8$. Thus, from Corollary \ref{corl-3} (Part (e)) we obtain $S_{89}=\left\{1-u_{89}^{-}, u_{89}^{+}\right\}$ and $I_{89}=\left[1-u_{89}^{-}, u_{89}^{+}\right]$, where $u_{89}^{-}$ and $u_{89}^{+}$ are calculated by $u_{89}^{-}=1-\sqrt{\left(1-b_{8}\right)^{2}-\left(1-a_{89}^{-}\right)^{2}}=1$ and $u_{89}^{+}=1-\sqrt{\left(1-b_{8}\right)^{2}-\left(1-a_{89}^{+}\right)^{2}}=0.8$ (see Table \ref{tbl_a3}). Therefore, $S_{89}=\{0,0.8\}$ and $I_{89}=[0,0.8]$. Also, see Remark \ref{rmk-2}, that indicates $u_{i j}^{+}=b_{i}$ if $0<b_{i}<a_{i j}^{+}=1$ and $u_{i j}^{-}=1$ if $a_{i j}^{-}=b_{i}$ for any arbitrary continuous Archimedean t-norm. By the similar calculations, it follows that $I_{i 9}=[0,1]$ for $i \in\{1, \ldots, 5\} \bigcup\{7,9,10\}$ and $I_{69}=[0.2,0.8]$. So, from Definition \ref{def-4}, we have $I_{9}=\bigcap_{i \in \mathscr{I}} I_{i 9}=[0.2,0.8]$ and $S_{89}^{\prime}=S_{89} \cap I_{9}=\{0.8\}$. In Section 6, Tables $1-4$ show all sets $I_{i j}$, $S_{i j}, I_{j}$ and $S_{i j}^{\prime}$ for each $i \in \mathscr{I}$ and $j \in \mathscr{J}$, respectively.
\end{example}
\section{Resolution of the feasible solution set for Problem (1)}\label{sec-3}
The following lemma gives two necessary conditions for the feasibility of Problem (\ref{eq_1}).

\begin{lemma}\label{lm-4}
	Suppose that $\varphi$ is a continuous Archimedean t-norm.
	\begin{itemize}
		\item[\textbf{(a)}] If $S\left(A^{+}, A^{-}, b\right) \neq \varnothing$, then $I_{j} \neq \varnothing$, $\forall j \in \mathscr{J}$.
		\item[\textbf{(b)}] If $S\left(A^{+}, A^{-}, b\right) \neq \varnothing$, then for each $i \in \mathscr{I}$, there exists at least one $j_{i} \in \mathscr{J}$ such that $S_{i j_{i}}^{\prime} \neq \varnothing$.
	\end{itemize}
\end{lemma}

\begin{proof}
	\textbf{(a)} By contradiction, suppose that $S\left(A^{+}, A^{-}, b\right) \neq \varnothing$ and $I_{j_{0}}=\varnothing$ for some $j_{0} \in \mathscr{J}$. Hence, there exists $i_{0} \in \mathscr{I}$ such that $x \notin I_{i_{0} j_{0}}$ (Definition \ref{def-4}), $\forall x \in S$ $\left(A^{+}, A^{-}, b\right)$. So, from Definition \ref{def-3}, we have $$\max \left\{\varphi\left(a_{i_{0} j_{0}}^{+}, x_{j_{0}}\right), \varphi\left(a_{i_{0} j_{0}}^{-}, 1-x_{j_{0}}\right)\right\}>b_{i_{0}}$$ that implies $x \notin S_{i_{0}}$. Consequently, since $S\left(A^{+}, A^{-}, b\right)=\bigcap_{i \in \mathscr{I}} S_{i}$, we obtain $x \notin S\left(A^{+}, A^{-}, b\right)$ that is a contradiction. \textbf{(b)} By contradiction, assume that $x \in S$ $\left(A^{+}, A^{-}, b\right)$ and there exists some $i_{0} \in \mathscr{I}$ such that $S_{i_{0} j}^{\prime}=\varnothing$, $\forall j \in \mathscr{J}$. So, $x_{j} \notin S_{i_{0} j}^{\prime}=S_{i_{0} j} \cap I_{j}$, $\forall j \in \mathscr{J}$ (Definition \ref{def-4}). Now, if $x_{j} \notin I_{j}$ for some $j \in \mathscr{J}$, then from Part (a) we obtain $S\left(A^{+}, A^{-}, b\right) \neq \varnothing$, that is a contradiction. On the other hand, if $x_{j} \notin S_{i_{0} j}^{\prime}$ and $x_{j} \in I_{j}$, $\forall j \in \mathscr{J}$, then it is concluded that $x_{j} \notin S_{i_{0} j}$, $\forall j \in \mathscr{J}$, that implies $\max \left\{\varphi\left(a_{i_{0}j}^{+}, x_{j}\right), \varphi\left(a_{i_{0}j}^{-}, 1-x_{j}\right)\right\}<b_{i_{0}}$, $\forall j \in \mathscr{J}$ (Definition \ref{def-3}). Therefore, $x \notin S_{i_{0}}$, that contradicts the assumption that $x \in S$ $\left(A^{+}, A^{-}, b\right)$.
\end{proof}

The following lemma provides a necessary and sufficient condition to guarantee that a given $x \in[0,1]^{n}$ is feasible for Problem (\ref{eq_1}) or not.

\begin{lemma}\label{lm-5}
	Suppose that $\varphi$ is a continuous Archimedean t-norm. Then, $x \in S\left(A^{+}, A^{-}, b\right)$ if and only if the following statements hold true:
	\begin{itemize}
		\item[\textbf{(I)}]
		$x_{j} \in I_{j}$, $\forall j \in \mathscr{J}$.
		\item[\textbf{(II)}]
		For each $i \in \mathscr{I}$, there exists at least one $j_{i} \in \mathscr{J}$ such that $x_{j_{i}} \in S_{i j_{i}}^{\prime}$.
	\end{itemize}
\end{lemma}

\begin{proof}
	Suppose that $x \in[0,1]^{n}$ satisfies the conditions (I) and (II). So, from the condition (I) and Definitions 3 and 4, we have $x_{j} \in I_{i j}$ and therefore $\max \left\{\varphi\left(a_{i j}^{+}, x_{j}\right), \varphi\left(a_{i j}^{-}, 1-x_{j}\right)\right\} \leq b_{i}$, $\forall i \in \mathscr{I}$ and $\forall j \in \mathscr{J}$. On the other hand, from the condition (II), and Definitions 3 and 4, for each $i \in \mathscr{I}$ we have $\max \left\{\varphi\left(a_{ij_{i}}^{+}, x_{j_{i}}\right), \varphi\left(a_{ij_{i}}^{-}, 1-x_{j_{i}}\right)\right\}=b_{i}$ for some $j_{i} \in \mathscr{J}$. Consequently, for each $i \in \mathscr{I}$, it is	concluded that $$\max _{j=1}^{n}\left\{\max \left\{\varphi\left(a_{i j}^{+}, x_{j}\right), \varphi\left(a_{i j}^{-}, 1-x_{j}\right)\right\}\right\}=\max \left\{\varphi\left(a_{i j_{i}}^{+}, x_{j_{i}}\right), \varphi\left(a_{i j_{i}}^{-}, 1-x_{j_{i}}\right)\right\}=b_{i}.$$ Thus, $x \in S_{i}$, $\forall i \in \mathscr{I}$, that implies $x \in \bigcap_{i \in \mathscr{I}} S_{i}=S\left(A^{+}, A^{-}, b\right)$. The converse statement is obtained by reversing the argument.
\end{proof}

\begin{remark}\label{rmk-4}
	From Lemma \ref{lm-5}, it follows that if $x \in S\left(A^{+}, A^{-}, b\right)$, then $L_{j} \leq x_{j} \leq U_{j}$, $\forall j \in \mathscr{J}$ , where $L_{j}$ and $U_{j}$ are the lower and upper bound of $I_{j}=\left[L_{j}, U_{j}\right]$, respectively.
\end{remark}

\begin{corollary}\label{corl-5}
	For each $i \in \mathscr{I}, x \in S_{i}$ if and only $x_{j} \in I_{i j}$, $\forall j \in \mathscr{J}$, and there exists at least one $j_{i} \in \mathscr{J}$ such that $x_{j_{i}} \in S_{i j_{i}}^{\prime}$.
\end{corollary}

\begin{proof}
	The proof is resulted from Lemma \ref{lm-5}, where $S\left(A^{+}, A^{-}, b\right)$ and $I_{j}$ are replaced by $S_{i}$ and $I_{i j}$, respectively.
\end{proof}

\begin{definition}\label{def-5}
	For each $i \in \mathscr{I}$, we define $\mathscr{J}_{i}=\left\{j \in \mathscr{J}: S_{i j}^{\prime} \neq \varnothing\right\}$. Similarly, for each $j \in \mathscr{J}$ , let $\mathscr{I}_{j}=\left\{i \in \mathscr{I}: S_{i j}^{\prime} \neq \varnothing\right\}$.
\end{definition}

\begin{definition}\label{def-6}
	A function $e$ (on $\mathscr{I}$ ) is said to be an admissible function if $e(i) \in \mathscr{J}_{i}(e)$, $\forall i \in \mathscr{I}$, where
	\begin{itemize}
		\item[\textbf{(I)}]
		$\mathscr{J}_{1}(e)=\mathscr{J}_{1}$.
		\item[\textbf{(II)}]
		$\mathscr{I}_{j}(e, i)=\{k \in \mathscr{I}: 1 \leq k<i ~ \text{and}~ e(k)=j\}$, $\forall j \in \mathscr{J}$ and $\forall i \in \mathscr{I}-\{1\}$.
		\item[\textbf{(III)}]
		$\mathscr{J}_{i}(e)=\left\{j \in \mathscr{J}_{i}: \mathscr{I}_{j}(e, i)=\varnothing\right. ~ \text{or}  \left.S_{i j}^{\prime} \cap\left(\bigcap_{k \in \mathscr{I}_{j}(e, i)} S_{k j}^{\prime}\right) \neq \varnothing\right\}$, $\forall i \in \mathscr{I}-\{1\}$.
	\end{itemize}
	Also, let $E$ be the set of all the admissible functions. For the sake of convenience, we can represent each $e$ as the vector $e=\left[j_{1}, \ldots, j_{m}\right]$ in which $e(i)=j_{i}$, $\forall i \in \mathscr{I}$.
\end{definition}

\begin{definition}\label{def-7}
	For each $e \in E$, let $\mathscr{I}_{j}(e)=\{i \in \mathscr{I}: e(i)=j\}$ and $S(e)$ be the set of all the vectors $x=\left(x_{1}, \ldots, x_{n}\right)$ such that
	\begin{equation}\label{eq-11}
		x_{j} \in\left\{\begin{array}{ll}
			\bigcap_{i \in \mathscr{I}_{j}(e)} S_{i j}^{\prime} & , \mathscr{I}_{j}(e) \neq \varnothing \\
			I_{j}                                               & , \mathscr{I}_{j}(e)=\varnothing
		\end{array} \quad, \forall j \in \mathscr{J}\right.
	\end{equation}
\end{definition}

\begin{corollary}\label{corl-6}
	Let $e: \mathscr{I} \rightarrow \bigcup_{i \in \mathscr{I}} \mathscr{J}_{i}$ be a function so that $e(i) \in \mathscr{J}_{i}$, $\forall i \in \mathscr{I}$. Then, $e \in E$ if and only if $\bigcap_{i \in \mathscr{I}_{j}(e)} S_{i j}^{\prime} \neq \varnothing$ for each $j \in \mathscr{J}$ such that $\mathscr{I}_{j}(e) \neq \varnothing$.
\end{corollary}

\begin{proof}
	See Appendix B.
\end{proof}

\begin{remark}\label{rmk-5}
	According to Corollary \ref{corl-6}, the number of admissible functions is bounded above by $\prod_{i \in \mathscr{I}}\left|\mathscr{J}_{i}\right|$, where $\left|\mathscr{J}_{i}\right|$ denotes the cardinality of $\mathscr{J}_{i}$. However, the actual number of admissible functions is usually much less than this value.
\end{remark}

Through the use of admissible functions, the following theorem determines the feasible solutions set for Problem (\ref{eq_1}).

\begin{theorem}\label{thm-3}
	Suppose that $\varphi$ is a continuous Archimedean t-norm. Then, $S\left(A^{+}, A^{-}, b\right)=\bigcup_{e \in E} S(e)$.
\end{theorem}

\begin{proof}
	See Appendix B.
\end{proof}

According to (\ref{eq-11}), for each $e \in E,~ x \in S(e)$ and $j \in \mathscr{J}$, we have either $x_{j} \in \bigcap_{i \in \mathscr{I}_{j}(e)} S_{i j}^{\prime}$ (if $\mathscr{I}_{j}(e) \neq \varnothing$ ) or $x_{j} \in I_{j}=\left[L_{j}, U_{j}\right]$ (if $\mathscr{I}_{j}(e)=\varnothing$ ). On the other hand, in the former case, we have $\bigcap_{i \in \mathscr{I}_{j}(e)} S_{i j}^{\prime} \neq \varnothing$ (from Corollary \ref{corl-6}) and also $\bigcap_{i \in \mathscr{I}_{j}(e)} S_{i j}^{\prime}$ is as one of the forms of $\left\{L_{j}\right\},~\left\{U_{j}\right\},~\left\{L_{j}, U_{j}\right\}$ or $\left[L_{j}, U_{j}\right]$ (from Corollary \ref{corl-4}). So, it is concluded that $S(e)$ ( $\forall e \in E$ ) is indeed a closed-bounded subset (and therefore, a compact subset) of $[0,1]^{n}$. This fact together with Theorem \ref{thm-3} implies that the feasible solution set of Problem (\ref{eq_1}) is determined by the union of a finite number of compact sets.

\begin{example}\label{ex-2}
	Consider the problem stated in Example \ref{ex-1}. According to Definition \ref{def-5} and Table \ref{t-4} (in Section 6), we obtain $\mathscr{J}_{1}=\{1,2,5\}$, $\mathscr{J}_{2}=\{1,8\}$, $\mathscr{J}_{3}=\{6,7,8\}$, $\mathscr{J}_{4}=\{2,3,4\}$, $\mathscr{J}_{5}=\{1,3,4,5\}$, $\mathscr{J}_{6}=\{8,9\}$, $\mathscr{J}_{7}=\{3,8,10\}$, $\mathscr{J}_{8}=\{9\}$, $\mathscr{J}_{9}=\{3,10\}$ and $\mathscr{J}_{10}=\{6,8\}$. Hence, according to Remark \ref{rmk-5}, the number of admissible functions is bounded above by $\prod_{i \in \mathscr{I}}\left|\mathscr{J}_{i}\right|=3 \times 2 \times 3 \times 3 \times 4 \times 2 \times 3 \times 1 \times 2 \times 2=5184$. Now, noting Definition \ref{def-6} and Corollary \ref{corl-6}, consider three functions $e_{1}=[5,1,6,3,5,8,10,9,10,6]$, $e_{2}=[1,8,6,2,3,8,10,9,10,6]$ and $e_{3}=[1,8,6,2,1,8,10,9,10,6]$ from $\mathscr{I}$ to $\bigcup_{i \in \mathscr{I}} \mathscr{J}_{i}$ so that $e_{p}(i) \in \mathscr{J}_{i}$ for each $i \in \mathscr{I}$ and each $p \in\{1,2,3\}$. For function $e_{1}$, we have $e_{1}(1)=e_{1}(5)=5$; $\mathscr{I}_{5}\left(e_{1}, 5\right)=\{1\}$ (from Part (II) of Definition \ref{def-6}); $S_{15}^{\prime}=\{0.6\}$ and $S_{55}^{\prime}=\{0.3\}$ (from Table \ref{t-4}). Therefore, $\mathscr{I}_{5}\left(e_{1}, 5\right) \neq \varnothing$ and $S_{55}^{\prime} \cap\left(\bigcap_{k \in \mathscr{I}_{5}\left(e_{1}, 5\right)} S_{k 5}^{\prime}\right)=S_{55}^{\prime} \cap S_{15}^{\prime}=\varnothing$, which implies $5 \notin \mathscr{J}_{5}\left(e_{1}\right)$ (Definition \ref{def-6}, Part (III)). So, since $e_{1}(5)=5 \notin \mathscr{J}_{5}\left(e_{1}\right)$, from Definition \ref{def-6} it follows that $e_{1}$ is not an admissible function. However, both $e_{2}$ and $e_{3}$ are admissible. For instance, for admissible function $e_{2}$ we have $\mathscr{I}_{1}\left(e_{2}\right)=\{1\}$, $\mathscr{I}_{2}\left(e_{2}\right)=\{4\}$, $\mathscr{I}_{3}\left(e_{2}\right)=\{5\}$, $\mathscr{I}_{4}\left(e_{2}\right)=\varnothing$, $\mathscr{J}_{5}\left(e_{2}\right)=\varnothing$, $\mathscr{I}_{6}\left(e_{2}\right)=\{3,10\}$, $\mathscr{I}_{7}\left(e_{2}\right)=\varnothing$, $\mathscr{I}_{8}\left(e_{2}\right)=\{2,6\}$, $\mathscr{I}_{9}\left(e_{2}\right)=\{8\}$ and $\mathscr{I}_{10}\left(e_{2}\right)=\{7,9\}$. So, by noting Table \ref{t-4}, it follows that $\bigcap_{i \in \mathscr{I}_{1}\left(e_{2}\right)} S_{i 1}^{\prime}=S_{11}^{\prime}=\{0.4\}$, $\bigcap_{i \in \mathscr{I}_{2}\left(e_{2}\right)} S_{i 2}^{\prime}=S_{42}^{\prime}=\{0.64\}$, $\bigcap_{i \in \mathscr{I}_{3}\left(e_{2}\right)} S_{i 3}^{\prime}=S_{53}^{\prime}=\{0.6\}$, $\bigcap_{i \in \mathscr{I}_{6}\left(e_{2}\right)} S_{i 6}^{\prime}=S_{36}^{\prime} \cap S_{10,6}^{\prime}=\{0,1\}$, $\bigcap_{i \in \mathscr{I}_{8}\left(e_{2}\right)} S_{i 8}^{\prime}=S_{28}^{\prime} \cap S_{68}^{\prime}=\{0.2\}$, $\bigcap_{i \in \mathscr{I}_{9}\left(e_{2}\right)} S_{i 9}^{\prime}=S_{89}^{\prime}=\{0.8\}$ and $\bigcap_{i \in \mathscr{I}_{10}\left(e_{2}\right)} S_{i, 10}^{\prime}=S_{7,10}^{\prime} \cap S_{9,10}^{\prime}=\{0.6\}$.
	
	Furthermore, from Table \ref{t_3}, $I_{4}=[0,1]$, $I_{5}=[0.3,0.6]$ and $I_{7}=[0,0.55]$. Hence, from Definition \ref{def-7}, $S\left(e_{2}\right)$ is calculated as the Cartesian product $S\left(e_{2}\right)=\{0.4\} \times\{0.64\} \times\{0.6\} \times[0,1] \times[0.3,0.6] \times\{0,1\} \times[0,0.55] \times\{0.2\} \times\{0.8\} \times\{0.6\}$.
\end{example}
\section{Optimal solutions of Problem (1)}\label{sec-4}
In this section, we introduce a subset of the feasible regions that contain optimal solutions. So, by considering only this new subset and ignoring other feasible solutions, we can accelerate the process of obtaining optimal solutions by reducing the search domain.

\begin{definition}\label{def-8}
	Suppose that $S\left(A^{+}, A^{-}, b\right) \neq \varnothing$. For each $e \in E$, let $x(e)=\left(x(e)_{1}, \ldots, x(e)_{n}\right)$ where for each $j \in \mathscr{J}$, component $x(e)_{j}$ is defined as follows:
	\begin{equation}\label{eq-12}
		x(e)_{j}= \begin{cases}\min \left\{\bigcap_{i \in \mathscr{I}_{j}(e)} S_{i j}^{\prime}\right\} & , \mathscr{I}_{j}(e) \neq \varnothing \\ L_{j} & , \mathscr{I}_{j}(e)=\varnothing\end{cases}
	\end{equation}
	where $L_{j}$ is the lower bound of $I_{j}=\left[L_{j}, U_{j}\right]$. Also, let $F=\{x(e): e \in E\}$.
\end{definition}

\begin{theorem}\label{thm-4}
	Suppose that $S\left(A^{+}, A^{-}, b\right) \neq \varnothing$ and $S^{*}$ denotes the set of optimal solutions for Problem (\ref{eq_1}). Then, $S^{*} \subseteq F \subseteq S\left(A^{+}, A^{-}, b\right)$.
\end{theorem}

\begin{proof}
	Clearly, $S^{*} \subseteq S\left(A^{+}, A^{-}, b\right)$. Also, from (\ref{eq-11}) and (\ref{eq-15}), it follows that if $x\left(e_{0}\right) \in F$ (for some $e_{0} \in E$ ), then $x\left(e_{0}\right) \in S\left(e_{0}\right) \subseteq \bigcup_{e \in E} S(e)$. So, by Theorem \ref{thm-3} , it is concluded that $F \subseteq S\left(A^{+}, A^{-}, b\right)$. To complete the proof, it is sufficient to show that any solution $x \in S\left(A^{+}, A^{-}, b\right)-F$ cannot be an optimal solution. For this purpose, let $x \in S\left(A^{+}, A^{-}, b\right)-F$. So, from Theorem \ref{thm-3}, there exists at least one $e_{0} \in E$ such that $x \in S\left(e_{0}\right)$. Now, by (\ref{eq-11}) and (\ref{eq-15}), we have $x\left(e_{0}\right)_{j} \leq x_{j}$, $\forall j \in \mathscr{J}$; and since $x \notin F$, at least one of the preceding inequalities must be strict. Without loss of generality, assume that $x\left(e_{0}\right)_{j_{0}}<x_{j_{0}}$. Hence, $c_{j_{0}} x\left(e_{0}\right)_{j_{0}} +\sum_{j \in \mathscr{J}-\left\{j_{0}\right\}} c_{j} x\left(e_{0}\right)_{j}<c_{j_{0}} x_{j_{0}}+\sum_{j \in \mathscr{J}-\left\{j_{0}\right\}} c_{j} x_{j}$, which implies $x \notin S^{*}$.
\end{proof}

Theorem \ref{thm-4} also provides a necessary optimality condition in the sense that if $x^{*}$ is an optimal solution to Problem (\ref{eq_1}), then it must belong to $F$. In other words, each solution $x(e) \in F$ is an optimal candidate solution.
\section{Simplification techniques}\label{sec-5}
In this section, by leveraging the special structure identified in the previous sections, some procedures are introduced that reduce the size of the original problem. Throughout this section, $S_{\{i\}}\left(A^{+}, A^{-}, b\right)$ refers to the feasible region of the reduced problem obtained by discarding the $i$'th equation from Problem (\ref{eq_1}). In the case of simpilicity, the proof of lemmas 6-13 have been presented in Appendix B.

5.1 Simplification techniques for finding the feasible solutions set

This subsection describes some simplification techniques to facilitate the determination of feasible regions. In the case of certain assumptions, all the techniques described in this subsection are derived from the fact that, any solution $x \in S_{\{i\}}\left(A^{+}, A^{-}, b\right)$ also satisfy the $i$'th equation of Problem (\ref{eq_1}). Therefore, in such a case, the $i$'th equation is apparent to be a redundant (irrelevant) constraint (i.e., a constraint that is not likely to affect the feasible region), and hence it can be omitted from further consideration.

\begin{lemma}\label{lm-6}
	Suppose that $S\left(A^{+}, A^{-}, b\right) \neq \varnothing$ and $i_{0} \in \mathscr{I}$. If $x \in S_{\left\{i_{0}\right\}}\left(A^{+}, A^{-}, b\right)$, then $x_{j} \in I_{j}$, $\forall j \in \mathscr{J}$. Particularly, $x_{j} \in I_{i_{0} j}$, $\forall j \in \mathscr{J}$.
\end{lemma}

\begin{lemma}\label{lm-7}
	Suppose that $S\left(A^{+}, A^{-}, b\right) \neq \varnothing$ and $i_{0} \in \mathscr{I}$. If $b_{i_{0}}=0$, then $S\left(A^{+}, A^{-}, b\right)=S_{\left\{i_{0}\right\}}\left(A^{+}, A^{-}, b\right)$.
\end{lemma}

\begin{corollary}\label{corl-7}
	Suppose that $S\left(A^{+}, A^{-}, b\right) \neq \varnothing$ and $i_{0} \in \mathscr{I}$. If $b_{i_{0}}=0$, then the $i_{0}$'th equation is a redundant constraint and it can be deleted.
\end{corollary}

\begin{lemma}\label{lm-8}
	Suppose that $S\left(A^{+}, A^{-}, b\right) \neq \varnothing$ and there exist $i_{0} \in \mathscr{I}$ and $j_{0} \in \mathscr{J}$ such that $I_{j_{0}}=\{k\}$ is a singleton set and $k \in S_{i_{0} j_{0}}^{\prime}$. Then, $x_{j_{0}}=k$, $\forall x \in S\left(A^{+}, A^{-}, b\right)$, and $S\left(A^{+}, A^{-}, b\right)=S_{\left\{i_{0}\right\}}\left(A^{+}, A^{-}, b\right)$.
\end{lemma}

Under the assumptions of Lemma \ref{lm-8}, the $i_{0}$'th equation is a redundant constraint. Lemma \ref{lm-8} additionally provides a further idea that is concerned with the case in which certain $x_{j}$'s can be calculated straight away without solving the problem, but just by reflecting on the special characteristics of the problem. Due to this, in such cases, any parts involved in these $x_{j}$'s can be excluded. These considerations have led us to come up with the following reduction rule.

\begin{corollary}\label{corl-8}
	Suppose that $S\left(A^{+}, A^{-}, b\right) \neq \varnothing$ and there exists $j_{0} \in \mathscr{J}$ such that $I_{j_{0}}=\{k\}$ is a singleton set. Then, $x_{j_{0}}=k$ for each feasible solution $x$. Also, the $j_{0}$'th column and any equation $i$ such that $k \in S_{ij_{0}}^{\prime}$ can be removed from the problem.
\end{corollary}

\begin{lemma}\label{lm-9}
	Suppose that $S\left(A^{+}, A^{-}, b\right) \neq \varnothing$ and there exist $i, i_{0} \in \mathscr{I}$ such that $S_{i j}^{\prime} \subseteq S_{i_{0} j}^{\prime}$, $\forall j \in \mathscr{J}$. Then, $S\left(A^{+}, A^{-}, b\right)=S_{\left\{i_{0}\right\}}\left(A^{+}, A^{-}, b\right)$.
\end{lemma}

\begin{corollary}\label{corl-9}
	Suppose that $S\left(A^{+}, A^{-}, b\right) \neq \varnothing$ and there exist $i, i_{0} \in \mathscr{I}$ such that $S_{i j}^{\prime} \subseteq S_{i_{0}j}^{\prime}$, $\forall j \in \mathscr{J}$. Then, the $i_{0}$'th equation is a redundant constraint and it can be deleted.
\end{corollary}

\begin{lemma}\label{lm-10}
	Suppose that $S\left(A^{+}, A^{-}, b\right) \neq \varnothing$ and there exist $i_{0} \in \mathscr{I}$ and $j_{0} \in \mathscr{J}$ such that $\mathscr{J}_{i_{0}}=\left\{j_{0}\right\}$ and $S_{i_{0} j_{0}}^{\prime}=\{k\}$ are singleton sets. Then, $x_{j_{0}}=k$, $\forall x \in S\left(A^{+}, A^{-}, b\right)$. Also, for each $i \in \mathscr{I}$ such that $k \in S_{i j_{0}}^{\prime}$, $S\left(A^{+}, A^{-}, b\right)=S_{\{i\}}\left(A^{+}, A^{-}, b\right) \cap\left\{x \in[0,1]^{n}: x_{j_{0}}=k\right\}$.
\end{lemma}

\begin{corollary}\label{corl-10}
	Suppose that $S\left(A^{+}, A^{-}, b\right) \neq \varnothing$ and there exist $i_{0} \in \mathscr{I}$ and $j_{0} \in \mathscr{J}$ such that $\mathscr{J}_{i_{0}}=\left\{j_{0}\right\}$ and $S_{i_{0} j_{0}}^{\prime}=\{k\}$ are singleton sets. Then, $x_{j_{0}}=k$ for each feasible solution $x$. Also, the $j_{0}$'th column and any equation $i$ such that $k \in S_{ij_{0}}^{\prime}$ can be removed from the problem.
\end{corollary}

5.2 Simplification techniques for finding optimal solutions

The techniques described in this subsection are based on Theorem \ref{thm-4}. These techniques reduce the size of the problem in such a way that either the set $F$ does not change or the search domain is reduced to a subset of $F$ that contains optimal solutions. In fact, by using these methods, some parts of the feasible region that do not contain optimal solutions may be removed. Therefore, these methods cannot be used to find the feasible region and are only suitable for finding optimal solutions.

\begin{lemma}\label{lm-11}
	Suppose that $S\left(A^{+}, A^{-}, b\right) \neq \varnothing$ and there exist $i_{0} \in \mathscr{I}$ and $j_{0} \in \mathscr{J}$ such that $\left|S_{i_{0} j_{0}}^{\prime}\right|=2$, where $\left|S_{i_{0} j_{0}}^{\prime}\right|$ denotes the cardinality of $S_{i_{0} j_{0}}^{\prime}$. If $x(e) \in S_{\left\{i_{0}\right\}}\left(A^{+}, A^{-}, b\right)$, then $x(e) \in S\left(A^{+}, A^{-}, b\right)$, $\forall e \in E$.
\end{lemma}

\begin{corollary}\label{corl-11}
	Suppose that $S\left(A^{+}, A^{-}, b\right) \neq \varnothing$ and there exist $i_{0} \in \mathscr{I}$ and $j_{0} \in \mathscr{J}$ such that $\left|S_{i_{0} j_{0}}^{\prime}\right|=2$. Then, the $i_{0}$'th equation is a redundant constraint and it can be deleted.
\end{corollary}

\begin{remark}\label{rmk-6}
	Let $i \in \mathscr{I}$ and $j \in \mathscr{J}$. After applying the above simplification technique, each non-empty set $S_{i j}$ is a singleton set. Consequently, if $\mathscr{I}^{\prime} \subseteq \mathscr{I}$ and $\bigcap_{i \in \mathscr{I}^{\prime}} S_{i j}^{\prime} \neq \varnothing$, then $S_{i_{1} j}^{\prime}=S_{i_{2} j}^{\prime}$ for each $i_{1} \in \mathscr{I}^{\prime}$ and each $i_{2} \in \mathscr{I}^{\prime}$.
\end{remark}

\begin{lemma}\label{lm-12}
	Suppose that $S\left(A^{+}, A^{-}, b\right) \neq \varnothing$ and there exists $j_{0} \in \mathscr{J}$ such that $L_{j_{0}} \in S_{ij_{0}}^{\prime}$, $\forall i \in \mathscr{I}_{j_{0}}$. Then, $x(e)_{j_{0}}=L_{j_{0}}$, $\forall e \in E$. Also, for each $i_{0} \in \mathscr{I}$ such that $L_{j_{0}} \in S_{i_{0} j_{0}}^{\prime}$, if $x(e) \in S_{\left\{i_{0}\right\}}\left(A^{+}, A^{-}, b\right)$, then $x(e) \in S\left(A^{+}, A^{-}, b\right)$, $\forall e \in E$.
\end{lemma}

\begin{corollary}\label{corl-12}
	Suppose that $S\left(A^{+}, A^{-}, b\right) \neq \varnothing$ and there exists $j_{0} \in \mathscr{J}$ such that $L_{j_{0}} \in S_{ij_{0}}^{\prime}$ , $\forall i \in \mathscr{I}_{j_{0}}$. If $x^{*}$ is an optimal solution of Problem (\ref{eq_1}), then $x_{j_{0}}^{*}=L_{j_{0}}$. Also, the $j_{0}$'th column and any equation $i$ such that $i \in \mathscr{I}_{j_{0}}$ can be removed from the problem.
\end{corollary}

\begin{proof}
	It is sufficient to prove $x_{j_{0}}^{*}=L_{j_{0}}$ for any optimal solution $x^{*}$. But, since $x^{*} \in S^{*}$, Theorem \ref{thm-4} implies that $x^{*} \in F$; that is, $x^{*}=x(e)$ for some $e \in E$. Now the result follows from Lemma \ref{lm-12}.
\end{proof}

\begin{lemma}\label{lm-13}
	Suppose that the simplification technique stated in Corollary \ref{corl-11} has been applied and there exist $j_{1}, j_{2} \in \mathscr{J}$ such that $\mathscr{I}_{j_{1}} \neq \varnothing$ and $\mathscr{I}_{j_{1}} \subseteq \mathscr{I}_{j_{2}}$. Also, define $F_{0}=\left\{x(e): e \in E_{0}\right\}$ and $E_{0}=\left\{e \in E: \mathscr{I}_{j_{1}}(e)=\varnothing\right\}$.
	\begin{itemize}
		\item[\textbf{(a)}] If $\bigcap_{i \in \mathscr{I}_{j_{2}}} S_{ij_{2}}^{\prime}=\left\{L_{j_{2}}\right\}$, then $S^{*} \subseteq F_{0} \subseteq F$.
		\item[\textbf{(b)}] If $\bigcap_{i \in \mathscr{I}_{j_1}} S_{ij_{1}}^{\prime}=\left\{U_{j_{1}}\right\}$, $\bigcap_{i \in \mathscr{I}_{j_{2}}} S_{i j_{2}}^{\prime}=\left\{U_{j_{2}}\right\}$ and $c_{j_{2}}\left(U_{j_{2}}-L_{j_{2}}\right)<c_{j_{1}}\left(U_{j_{1}}-L_{j_{1}}\right)$, then $S^{*} \subseteq F_{0} \subseteq F$.
	\end{itemize}
\end{lemma}

According to Lemma \ref{lm-13}, we can reduce our search to set $F_{0}$ instead of set $F$. To achieve this, we neglect all the admissible functions $e \in E-E_{0}$ by removing each $i \in \mathscr{I}_{j_{1}}$ before selecting the vectors $e$ to construct solutions $x(e)$. However, this is equivalent to assigning $x_{j_{2}}^{*}=L_{j_{2}}$ ($x_{j_{2}}^{*}=U_{j_{2}}$ in part (b) of Lemma \ref{lm-13}) and deleting the $j_{1}$'th column from the problem. The following reduction rule summarizes the previous discussion.

\begin{corollary}\label{corl-13}
	Suppose that the simplification technique stated in Corollary \ref{corl-11} has been applied.
	\begin{itemize} 
		\item[\textbf{(a)}] If there exist $j_{1}, j_{2} \in \mathscr{J}$ such that $\mathscr{I}_{j_{1}} \neq \varnothing$, $\mathscr{I}_{j_{1}} \subseteq \mathscr{I}_{j_{2}}$ and $\bigcap_{i \in \mathscr{I}_{j_{2}}} S_{i j_{2}}^{\prime}=\left\{L_{j_{2}}\right\}$, then $x_{j_{2}}^{*}=L_{j_{2}}$. Also, the $j_{1}$'th column can be removed from the problem.
		\item[\textbf{(b)}] If there exist $j_{1}, j_{2} \in \mathscr{J}$ such that $\mathscr{I}_{j_{1}} \subseteq \mathscr{I}_{j_{2}}$, $\bigcap_{i \in \mathscr{I}_{j_1}} S_{ij_{1}}^{\prime}=v \in\left\{L_{j_{1}}, U_{j_{1}}\right\}$, $\bigcap_{i \in \mathscr{I}_{j_{2}}} S_{ij_{2}}^{\prime}=\left\{U_{j_{2}}\right\}$ and $c_{j_{2}}\left(U_{j_{2}}-L_{j_{2}}\right)<c_{j_{1}}\left(v-L_{j_{1}}\right)$, then $x_{j_{2}}^{*}=U_{j_{2}}$. Also, the $j_{1}$'th column can be removed from the problem.
	\end{itemize}
\end{corollary}

\vskip1cm
5.3 Simplification strategies for finding optimal solutions

In this section, we describe two strategies that, like the techniques described in subsection 5.2, can only be used to find optimal solutions. Contrary to the previous simplification techniques, which were implemented at the beginning of the algorithm, these two strategies are used during the problem-solving process.

\begin{definition}\label{def-9}
	Suppose that the simplification technique stated in Corollary \ref{corl-11} has been applied. An admissible function $e$ is said to be a modified function if $e(i) \in \bar{\mathscr{J}}_{i}(e)$ , $\forall i \in \mathscr{I}$, where
	$$
	\bar{\mathscr{J}}_{i}(e)= \begin{cases}\min \left\{\mathscr{J}_{i}(e) \cap\{e(1), \ldots, e(i-1)\}\right\} & , \mathscr{J}_{i}(e) \cap\{e(1), \ldots, e(i-1)\} \neq \varnothing \\ \mathscr{J}_{i}(e) & , \text {otherwise}\end{cases}
	$$
	Also, let $\bar{E}$ be the set of all the modified functions and $\bar{F}=\{x(e): e \in \bar{E}\}$, where $x(e)$ is defined by (\ref{eq-15}).
\end{definition}

\begin{example}\label{ex-3}
	Consider admissible functions $e_{2}=[1,8,6,2,3,8,10,9,10,6]$ and $e_{3}=[1,8,6,2,1,8,10,9,10,6]$ stated in Example \ref{ex-2}. As calculated in Example \ref{ex-2}, $\mathscr{J}_{5}=\{1,3,4,5\}$. Also, for admissible function $e_{2}$, we have $e_{2}(1)=1$, $e_{2}(2)=8$, $e_{2}(3)=6$ and $e_{2}(4)=2$. Furthermore, Part (II) of Definition \ref{def-6} indicates that $\mathscr{I}_{1}\left(e_{2}, 5\right)=\{1\}$ and $\mathscr{I}_{3}\left(e_{2}, 5\right)=\mathscr{I}_{4}\left(e_{2}, 5\right)=\mathscr{I}_{5}\left(e_{2}, 5\right)=\varnothing$. On the other hand, from Table \ref{t-4}, it follows that $S_{51}^{\prime} \cap\left(\bigcap_{k \in \mathscr{I}_{1}\left(e_{2}, 5\right)} S_{k 1}^{\prime}\right)=S_{51}^{\prime} \cap S_{11}^{\prime}=\{0.4\}$. Therefore, from Part (III) of Definition \ref{def-6}, we obtain $\mathscr{J}_{5}\left(e_{2}\right)=\{1,3,4,5\}$, which implies $\mathscr{J}_{5}\left(e_{2}\right) \cap\left\{e_{2}(1), e_{2}(2), e_{2}(3), e_{2}(4)\right\}=\{1,3,4,5\} \cap\{1,8,6,2\}=\{1\} \neq \varnothing$. So, $\bar{\mathscr{J}}_{5}\left(e_{2}\right)=\{1\} \quad$ and $e_{2}(5)=3 \notin \bar{\mathscr{J}}_{5}\left(e_{2}\right)$, that means $e_{2}$ is not a modified function (Definition \ref{def-9}). However, function $e_{3}$ is an example of an admissible function that is also a modified function.
\end{example}

\begin{theorem}\label{thm-5}
	Suppose that the simplification technique stated in Corollary \ref{corl-11} has been applied. If $S\left(A^{+}, A^{-}, b\right) \neq \varnothing$ and $S^{*}$ denotes the set of optimal solutions for Problem (\ref{eq_1}), then $S^{*} \subseteq \bar{F} \subseteq F$.
\end{theorem}

\begin{proof}
	From Definition \ref{def-9}, it is clear that $\bar{F} \subseteq F$. To complete the proof, we shall show that for each $e \in E-\bar{E}$, there exists some $\bar{e} \in \bar{E}$ such that $\sum_{j \in \mathscr{J}} c_{j} x(\bar{e})_{j} \leq \sum_{j \in \mathscr{J}} c_{j} x(e)_{j}$. For this purpose, let $e \in E-\bar{E}$, $\mathscr{I}^{\prime}=\left\{i_{1}^{\prime}, \ldots, i_{k}^{\prime}\right\}$ and $\mathscr{J}^{\prime}=\left\{j_{1}^{\prime}, \ldots, j_{k}^{\prime}\right\}$ where $\mathscr{J}_{i_{p}^{\prime}}(e) \bigcap\left\{e(1), \ldots, e\left(i_{p}^{\prime}-1\right)\right\} \neq \varnothing$ and $e(i_{p}^{\prime})=j_{p}^{\prime} \in \mathscr{J}_{i_{p}^{\prime}}(e)-\left\{e(1), \ldots, e\left(i_{p}^{\prime}-1\right)\right\}$ for $p \in\{1, \ldots, k\}$. Now, consider the function $\bar{e}$ such that $\bar{e}(i)=\min \left\{\mathscr{J}_{i}(e) \cap\{e(1), \ldots, e(i-1)\}\right\}$, $\forall i \in \mathscr{I}^{\prime}$, and $\bar{e}(i)=e(i)$, otherwise. So, we have $\bar{e} \in \bar{E}$, $\mathscr{I}_{j}(\bar{e})=\mathscr{I}_{j}(e)$ $\left(\forall j \in \mathscr{J}-\mathscr{J}^{\prime}\right)$ and $\mathscr{I}_{j}(\bar{e}) \subseteq \mathscr{I}_{j}(e)$ $\left(\forall j \in \mathscr{J}^{\prime}\right)$. Therefore, $x(\bar{e})_{j}=x(e)_{j}$, $\forall j \in \mathscr{J}-\mathscr{J}^{\prime}$. Moreover, for each $j \in \mathscr{J}^{\prime}$ we have $x(e)_{j}=\min \left\{\bigcap_{i \in \mathscr{I}_{j}(e)} S_{i j}^{\prime}\right\} \in\left\{L_{j}, U_{j}\right\} \quad$ (from Corollary \ref{corl-4}), $x(\bar{e})_{j}=L_{j}$ if $\mathscr{I}_{j}(\bar{e})=\varnothing$ (from (\ref{eq-12})) and $x(\bar{e})_{j}=\min \left\{\bigcap_{i \in \mathscr{I}_{j}(\bar{e})} S_{i j}^{\prime}\right\}=\min \left\{\bigcap_{i \in \mathscr{I}_{j}(e)} S_{i j}^{\prime}\right\}=x(e)_{j}$ if $\mathscr{I}_{j}(\bar{e}) \neq \varnothing$ (from Remark \ref{rmk-6}). Hence, in any case, $x(\bar{e})_{j} \leq x(e)_{j}$ that implies $\sum_{j \in \mathscr{J}} c_{j} x(\bar{e})_{j} \leq \sum_{j \in \mathscr{J}} c_{j} x(e)_{j}$.
\end{proof}

Based on Theorem \ref{thm-5}, as the first simplification strategy, we can accelerate the identification of the optimal solutions set by focusing on the set $\bar{F}$ instead of $F$. In other words, those candidate solutions $x(e)$ (Theorem \ref{thm-4}) that belong to $F-\bar{F}$ are not optimal and therefore can be excluded from consideration. However, by using a branch-and-bound method (as the second simplification strategy), we can find an optimal solution without explicitly generating the whole set of solutions $x(e)$ $(e \in \bar{E})$. The branch-and-bound method incrementally builds modified functions $e \in \bar{E}$ by selecting the columns of $\bar{\mathscr{J}}_{i}(e)$ for each row $i$.
More precisely, let $n_{i}$ denote the cardinality of $\bar{\mathscr{J}}_{i}(e)$. The branch-and-bound method begins by choosing the columns of $\bar{\mathscr{J}}_{1}(e)$ and generating all partial sequences $e=\left[j_{1}\right]$, $\forall j_{1} \in \bar{\mathscr{J}}_{1}(e)$. Each $e$ is represented by one node, and therefore this yields $n_{1}$ partial branches denoted by nodes 1 to $n_{1}$. Also, for each $e=\left[j_{1}\right]$, we have $\mathscr{I}_{j_{1}}(e)=\{1\}$ and $\mathscr{I}_{j}(e)=\varnothing$, $\forall j \in \mathscr{J}-\left\{j_{1}\right\}$. So, associated with each node $e=\left[j_{1}\right]$, the solution $x(e)$ can be obtained by (\ref{eq-12}) with the objective value of $c x(e)$. Subsequently, we branch at the node with a smaller objective value. Without loss of generality, suppose that the solution $x\left(e^{*}\right)$, generated by $e^{*}=\left[j_{1}^{*}\right]$, has the smallest objective value. Then, in the second step, we focus on partial sequences $e=\left[j_{1}^{*}, j_{2}\right]$, $\forall j_{2} \in \bar{\mathscr{J}}_{2}(e)$, that yields nodes $n_{1}+1$ to $n_{1}+n_{2}$. Similarly, for new nodes $e=\left[j_{1}^{*}, j_{2}\right]$ $\left(j_{2} \in \bar{\mathscr{J}}_{2}(e)\right)$, solutions $x(e)$ are attained by (\ref{eq-12}), and we branch on the new node with a smaller objective value. It should be noted that for successive nodes $e^{\prime}=\left[j_{1}, \ldots, j_{p}\right]$ and $e^{\prime \prime}=\left[j_{1}, \ldots, j_{p}, j_{p+1}\right]$, either $j_{p+1} \in\left\{j_{1}, j_{2}, \ldots, j_{p}\right\}$ or $j_{p+1} \notin\left\{j_{1}, j_{2}, \ldots, j_{p}\right\}$. In the former case, we have $\mathscr{I}_{j}\left(e^{\prime}\right)=\mathscr{I}_{j}\left(e^{\prime \prime}\right)$, $\forall j \in \mathscr{J}$, which implies $x\left(e^{\prime}\right)=x\left(e^{\prime \prime}\right)$ (Definition \ref{def-8}). But, in the latter case, $\mathscr{I}_{j_{p+1}}\left(e^{\prime}\right)=\varnothing$, $\mathscr{I}_{j_{p+1}}\left(e^{\prime \prime}\right)=\{p+1\}$ and $\mathscr{I}_{j}\left(e^{\prime}\right)=\mathscr{I}_{j}\left(e^{\prime \prime}\right)$, $\forall j \in \mathscr{J}-\left\{j_{p+1}\right\}$, which together with Definition \ref{def-8} and Corollary \ref{corl-4} imply $x\left(e^{\prime}\right)_{j_{p+1}}=L_{j_{p+1}}$, $x\left(e^{\prime \prime}\right)_{j_{p+1}} \in\left\{L_{j_{p+1}}, U_{j_{p+1}}\right\}$ and $x\left(e^{\prime}\right)_{j}=x\left(e^{\prime \prime}\right)_{j}$, $\forall j \in \mathscr{J}-\left\{j_{p+1}\right\}$; that is, $x\left(e^{\prime}\right)_{j} \leq x\left(e^{\prime \prime}\right)_{j}$, $\forall j \in \mathscr{J}$. So, in any case, we have $c x\left(e^{\prime}\right) \leq c x\left(e^{\prime \prime}\right)$, which means the objective function is non-decreasing on any branches. By the same reasoning as before, we finally attain a node that corresponds to a complete sequence $e_{1}=\left[j_{1}^{*}, \ldots, j_{m-1}^{*}, j_{m}\right]$ in which $j_{m} \in \bar{\mathscr{J}}_{m}(e)$. Hence, $e_{1} \in \bar{E}$ and then $x\left(e_{1}\right) \in \bar{F}$. The resulting node can be viewed as a candidate solution with $z_{1}=c x\left(e_{1}\right)$. Now, since the objective function is non-decreasing, all nodes whose values cannot be lower than $z_{1}$, can be eliminated from consideration. However, we branch on those nodes so that it is still possible for them to yield a complete sequence $e_{2}$ such that $z_{1}>z_{2}$. Therefore, evaluating some nodes could result in the elimination of many possible solutions or branches from further consideration. Hence, if the solution of the partial branch we are evaluating is larger than the best solution obtained from a complete branch, there is no need to continue evaluating the current branch. The preceding discussion is now summarized as an algorithm.

\begin{algorithm}[!ht]
	\DontPrintSemicolon
	\SetKwInOut{Input}{Input}
	\SetKwInOut{Output}{Output}
	
	\Input{$A^{+}$, $A^{-} Metrices$ and Given continuous Archimedean t-norm $\varphi$}
	\Output{Find $x^{*}$ (An Optimal Solution)}
	
	Compute sets $I_{i j}, S_{i j}, I_{j}$ and $S_{i j}^{\prime}$ for each $i \in \mathscr{I}$ and each $j \in \mathscr{J}$ (Remark \ref{rmk-2} and Table \ref{tbl_a3}).
	
	\If{$I_{j}=\varnothing$ for some $j \in \mathscr{J}$,}{
		stop; the problem is infeasible (Lemma \ref{lm-4} (a)).
	}
	\Else{
		\If{$S_{i j}^{\prime}=\varnothing$ for some $i \in \mathscr{I}$ and each $j \in \mathscr{J}$,}{
			Stop; The Problem is Infeasible (Lemma \ref{lm-4} (b)).
		}
		\Else{
			Apply simplification techniques (Corollaries \ref{corl-7}-\ref{corl-13}).
			
			Delete the corresponding columns and redundant equations from the problem.
			
			Employ the branch-and-bound method on the modified functions $e \in \bar{E}$.
			
			Generate an optimal solution for the problem and determine the optimal value.
		}
	}
	\caption{Find Optimal Solution for Problem (\ref{eq_1})}
	\label{myalgorithm}
\end{algorithm}

\begin{figure}[!ht]
	\begin{center}
		\includegraphics[scale=0.5]{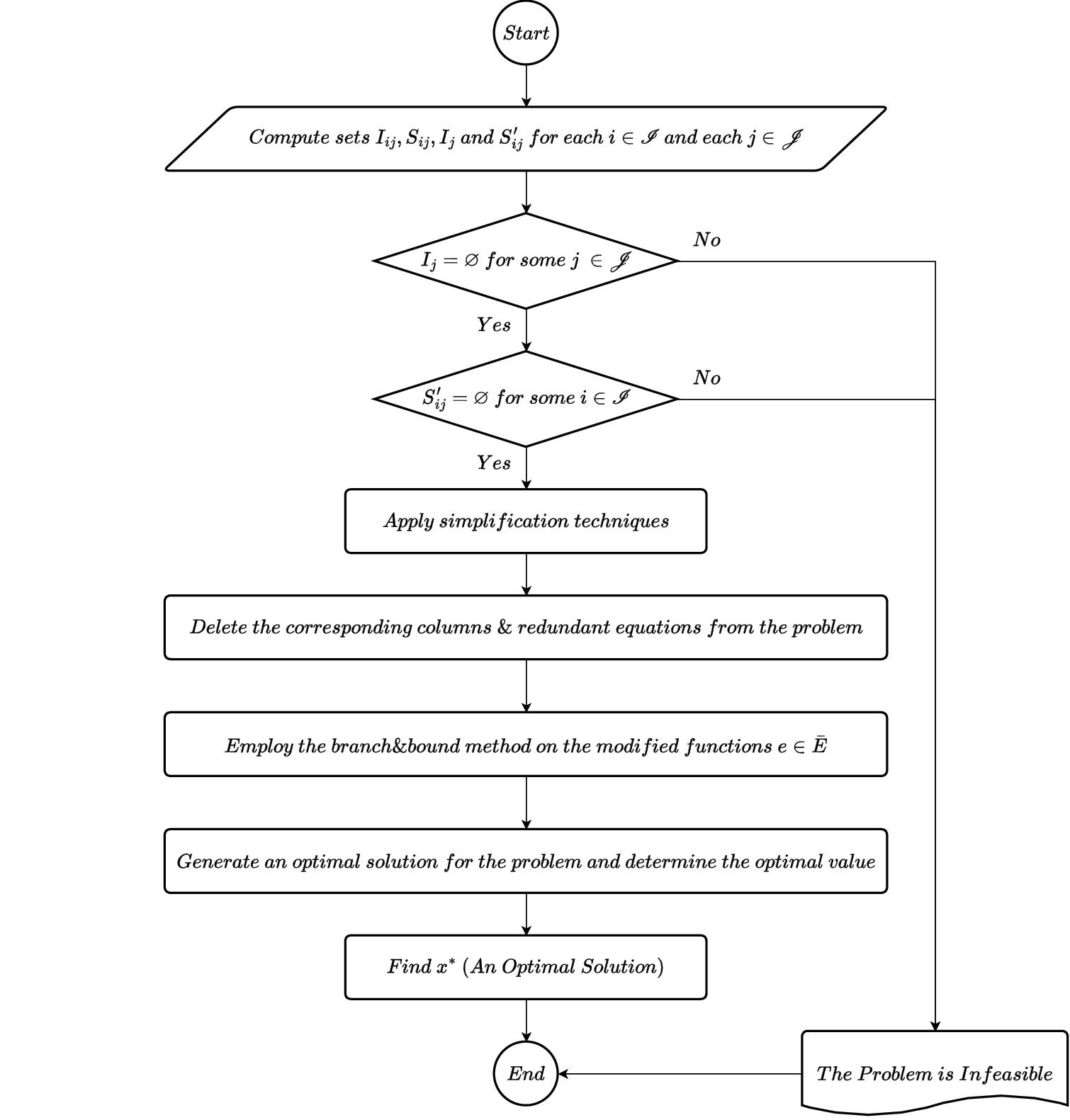}
		\captionsetup{justification=centering}
		\caption{Flowchart of Algorithm \ref{myalgorithm}}
		\label{fig_2}
	\end{center}
\end{figure}

In the first step of the algorithm, the $I_{i j}, S_{i j}, I_{j}$ and $S_{i j}^{\prime}$ sets are computed using at most $O(mn)$ comparisons with the aid of Remark \ref{rmk-2} and Table \ref{tbl_a3}. In addition, the conditions in steps 2 and 6 of algorithm involve a comparison of at most $O(mn)$. It may be necessary to perform up to $O(m(m+1)/2)$ operations in order to apply the simplification techniques in step 10 of algorithm. There is, however, a possibility that the simplification process may exhibit non-polynomial complexity if the problem variables exceed a polynomial growth rate. Similarly, step 11 of algorithm requires up to $O(mn)$ operations (considering that $m$ and $n$ have been modified as a result of the simplification techniques). In step 12 of algorithm, if the upper bound of the admissible function $|E|$ exceeds polynomial complexity, as mentioned in the introduction, the problem falls into the NP-hard complexity class of computation. In addition, Flowchart \ref{fig_2} provides a better understanding of Algorithm \ref{myalgorithm}. 

According to the algorithm described in this paper, we are able to obtain an exact (non-numerical and non-probabilistic) solution to Problem \ref{eq_1}. Unlike other studies that are focused on nonlinear objective functions, this paper examines a linear problem with constraints defined by bipolar fuzzy equations using continuous Archimedean t-norms. It is therefore not possible to analyze algorithms using benchmarks in the same manner as numerical or heuristic methods. Moreover, the study proved theorems, lemmas, and results common to all continuous Archimedean t-norms relevant to Problem \ref{eq_1} (Table \ref{tbl_a1} contains a list of common continuous Archimedean t-norms). In consequence, variation in the t-norms does not significantly alter the core steps of Algorithm \ref{myalgorithm}; rather, the differences lie in their properties and characteristics. According to the findings of the current study, future work may focus on enhancing the proposed algorithm or further simplifying Problem \ref{eq_1}, which would collectively contribute to the improved practical performance of the algorithm.

\section{Numerical example and application}\label{sec-6}
Consider a supplier of various products, identified as $p_{1}, p_{2}, \cdots, p_{n}$. The supplier seeks to optimize the public awareness and thereby attributes to all the products a degree of appreciation, $x_{1}, x_{2}, \cdots, x_{n}$. For the $j$'th product, the degree of appreciation $x_{j}$ is expressed as a real number in the unit interval $[0,1]$. Hence, its degree of disappreciation is equated to $1-x_{j}$. The degree of appreciation thus adopts a bipolar character on $[0,1]$, which means that value 0.5 means neutrality with regard to appreciation. Since all the products do not have the same importance, we assign the weight $c_{j}$ to the degree of appreciation $x_{j}$, and therefore the supplier's target can be considered as the optimization of the objective function $\sum_{j=1}^{n} c_{j} x_{j}$. Regarding the publicity effects of products, there is a basic threshold for appreciation (or disappreciation). As an example, assume that the products are promoted in $m$ areas/markets, denoted by $A_{1}, A_{2}, \cdots, A_{m}$. The threshold degree of appreciation of $p_{j}$ at the $i$'th market $A_{i}$ is $a_{i j}\left(a_{i j} \in[0,1]\right)$. If the degree of appreciation $x_{j}$ is less than $a_{i j}$, in that case, the people might not pay attention to the product. Hence we assume that people do not notice the product and the corresponding public awareness could be neglected when the degree of appreciation is below the threshold. If $x_{j}$ is no less than $a_{i j}$, then the effective public awareness is $x_{j}-a_{i j}$. As a result, the effective public awareness could be represented as $\max \left\{x_{j}-a_{i j}, 0\right\}$, which is the relative degree of appreciation of $p_{j}$ at the market $A_{i}$. Furthermore, if the threshold degree of disappreciation of $p_{j}$ at $A_{i}$ is $\tilde{a}_{i j}\left(\tilde{a}_{i j} \in[0,1]\right)$, then the relative degree of disappreciation could be represented as $\max \left\{\left(1-x_{j}\right)-\tilde{a}_{i j}, 0\right\}$. Let $a_{i j}^{+}=1-a_{i j} \in[0,1]$ and $a_{i j}^{-}=1-$ $\tilde{a}_{i j} \in[0,1]$. Then the relative degrees of appreciation and disappreciation could be written as $\max \left\{a_{i j}^{+}+x_{j}-1,0\right\}=\varphi\left(a_{i j}^{+}, x_{j}\right)$ and $\max \left\{a_{i j}^{-}+\left(1-x_{j}\right)-1,0\right\}=$ $\varphi\left(a_{i j}^{-}, 1-x_{j}\right)$, respectively, where $\varphi$ is the Lukasiewicz t-norm. Therefore, public awareness of product $p_{j}$ at market $A_{i}$ could be represented as $\max \left\{\varphi\left(a_{i j}^{+}, x_{j}\right), \varphi\left(a_{i j}^{-}, 1-x_{j}\right)\right\}$. Moreover, assume that market research has revealed that public awareness at a level of $b_{i}$ is the right value for the products to be sold at market $A_{i}$. Thus, achieving such a level of public awareness will result in increased sales for the company. For an arbitrary market $A_{i}$, there is at least one product, of which public awareness reaches the highest value $b_{i}$. In this regard, the above conditions could be expressed as $\max _{j \in J}\left\{\max \left\{\varphi\left(a_{i j}^{+}, x_{j}\right), \varphi\left(a_{i j}^{-}, 1-x_{j}\right)\right\}\right\}=$ $b_{i}, i=1,2, \cdots, m$. Hence, the above optimization model is reduced to Problem (\ref{eq_1}) in which $\varphi$ is the Lukasiewicz t- norm.

\begin{example}\label{ex-4}
	Here we will take a closer look at the problem presented in Examples \ref{ex-1}, \ref{ex-2}, and \ref{ex-3}.
	\begin{step}\label{st-1}
		Tables 1 and 2 show all sets $I_{i j}$ and $S_{i j}$ for each $i \in \mathscr{I}$ and $j \in \mathscr{J}$, respectively. In both tables, row $i ~(i \in \mathscr{I})$ corresponds to equation $i$ and column $j ~(j \in \mathscr{J})$ corresponds to variable $x_{j}$. By Table \ref{t_1} and Definition \ref{def-4}, all sets $I_{j}$ are shown in Table \ref{t_3} for each $j \in \mathscr{J}$. Also, all set $S_{i j}^{\prime}$ ($\forall i \in \mathscr{I}$ and $\forall j \in \mathscr{J})$ are summarized in Table \ref{t-4} by using Table \ref{t_2}, Table \ref{t_3} and Definition \ref{def-4}.
		
		\begin{table}[!ht]
			\centering
			\caption{Sets $I_{i j}=\left[L_{i j}, U_{i j}\right]$ for each $i \in \mathscr{I}$ and $j \in \mathscr{J}$.}
			\label{t_1}
			\begin{tabular}{|c|c|c|c|c|c|c|c|c|c|}
				\hline
				$[0.4,1]$ & $[0.4,1]$  & $[0,1]$   & $[0,1]$ & $[0,0.6]$ & $[0,1]$ & $[0,1]$    & $[0,1]$   & $[0,1]$     & $[0,1]$   \\
				\hline
				$[0,1]$   & $[0,1]$    & $[0,1]$   & $[0,1]$ & $[0,1]$   & $[0,1]$ & $[0,1]$    & $[0.2,1]$ & $[0,1]$     & $[0,1]$   \\
				\hline
				$[0,1]$   & $[0,1]$    & $[0,1]$   & $[0,1]$ & $[0,1]$   & $[0,1]$ & $[0,0.55]$ & $[0,1]$   & $[0,1]$     & $[0,1]$   \\
				\hline
				$[0,1]$   & $[0,0.64]$ & $[0,0.6]$ & $[0,1]$ & $[0,1]$   & $[0,1]$ & $[0,1]$    & $[0,1]$   & $[0,1]$     & $[0,1]$   \\
				\hline
				$[0.4,1]$ & $[0,1]$    & $[0,0.6]$ & $[0,1]$ & $[0.3,1]$ & $[0,1]$ & $[0,1]$    & $[0,1]$   & $[0,1]$     & $[0,1]$   \\
				\hline
				$[0,1]$   & $[0,1]$    & $[0,1]$   & $[0,1]$ & $[0,1]$   & $[0,1]$ & $[0,1]$    & $[0.2,1]$ & $[0.2,0.8]$ & $[0,1]$   \\
				\hline
				$[0,1]$   & $[0,1]$    & $[0,1]$   & $[0,1]$ & $[0,1]$   & $[0,1]$ & $[0,1]$    & $[0,1]$   & $[0,1]$     & $[0,0.6]$ \\
				\hline
				$[0,1]$   & $[0,1]$    & $[0,1]$   & $[0,1]$ & $[0,1]$   & $[0,1]$ & $[0,1]$    & $[0,1]$   & $[0,0.8]$   & $[0,1]$   \\
				\hline
				$[0,1]$   & $[0,1]$    & $[0,1]$   & $[0,1]$ & $[0,1]$   & $[0,1]$ & $[0,1]$    & $[0,1]$   & $[0,1]$     & $[0.6,1]$ \\
				\hline
				$[0,1]$   & $[0,1]$    & $[0,1]$   & $[0,1]$ & $[0,1]$   & $[0,1]$ & $[0,1]$    & $[0,1]$   & $[0,1]$     & $[0,1]$   \\
				\hline
			\end{tabular}%
		\end{table}
		
		\begin{table}[!ht]
			\centering
			\caption{Sets $S_{i j}$ for each $i \in \mathscr{I}$ and $j \in \mathscr{J}$.}
			\label{t_2}
			\begin{tabular}{|c|c|c|c|c|c|c|c|c|c|}
				\hline
				$\{0.4\}$     & $\{0.4\}$     & $\varnothing$ & $\varnothing$ & $\{0.6\}$     & $\varnothing$ & $\varnothing$ & $\{0\}$       & $\varnothing$ & $\{1\}$       \\
				\hline
				$\{1\}$       & $\varnothing$ & $\varnothing$ & $\varnothing$ & $\{0,1\}$     & $\varnothing$ & $\{1\}$       & $\{0.2\}$     & $\varnothing$ & $\varnothing$ \\
				\hline
				$\varnothing$ & $\varnothing$ & $\varnothing$ & $\varnothing$ & $\varnothing$ & $\{0,1\}$     & $\{0.55\}$    & $\{1\}$       & $\varnothing$ & $\varnothing$ \\
				\hline
				$\varnothing$ & $\{0.64\}$    & $\{0.6\}$     & $\{1\}$       & $\varnothing$ & $\varnothing$ & $\varnothing$ & $\varnothing$ & $\{0,1\}$     & $\{1\}$       \\
				\hline
				$\{0.4\}$     & $\varnothing$ & $\{0.6\}$     & $\{1\}$       & $\{0.3\}$     & $\varnothing$ & $\varnothing$ & $\{0\}$       & $\{0\}$       & $\{0\}$       \\
				\hline
				$\varnothing$ & $\varnothing$ & $\varnothing$ & $\varnothing$ & $\varnothing$ & $\varnothing$ & $\varnothing$ & $\{0.2\}$     & $\{0.2,0.8\}$ & $\varnothing$ \\
				\hline
				$\varnothing$ & $\{0\}$       & $\{0,1\}$     & $\varnothing$ & $\varnothing$ & $\varnothing$ & $\{1\}$       & $\{1\}$       & $\varnothing$ & $\{0.6\}$     \\
				\hline
				$\varnothing$ & $\varnothing$ & $\varnothing$ & $\varnothing$ & $\varnothing$ & $\varnothing$ & $\varnothing$ & $\varnothing$ & $\{0,0.8\}$   & $\varnothing$ \\
				\hline
				$\varnothing$ & $\varnothing$ & $\{0,1\}$     & $\varnothing$ & $\{0\}$       & $\varnothing$ & $\varnothing$ & $\{0\}$       & $\{1\}$       & $\{0.6\}$     \\
				\hline
				$\varnothing$ & $\varnothing$ & $\varnothing$ & $\varnothing$ & $\varnothing$ & $\{0,1\}$     & $\varnothing$ & $\{1\}$       & $\varnothing$ & $\varnothing$ \\
				\hline
			\end{tabular}%
		\end{table}
		
		\begin{table}[!ht]
			\centering
			\caption{Sets $I_{j}=\left[L_{j}, U_{j}\right]$ for each $j \in \mathscr{J}$.}
			\label{t_3}
			\begin{tabular}{|l|l|l|l|l|l|l|l|l|l|}
				\hline
				$[0.4,1]$ & $[0.4,0.64]$ & $[0,0.6]$ & $[0,1]$ & $[0.3,0.6]$ & $[0,1]$ & $[0,0.55]$ & $[0.2,1]$ & $[0.2,0.8]$ & $\{0.6\}$ \\
				\hline
			\end{tabular}%
		\end{table}
		
		\begin{table}[!ht]
			\centering
			\caption{Sets $S_{i j}^{\prime}$ for each $i \in \mathscr{I}$ and $j \in \mathscr{J}$.}
			\label{t-4}
			\begin{tabular}{|c|c|c|c|c|c|c|c|c|c|}
				\hline
				$\{0.4\}$     & $\{0.4\}$     & $\varnothing$ & $\varnothing$ & $\{0.6\}$     & $\varnothing$ & $\varnothing$ & $\varnothing$ & $\varnothing$ & $\varnothing$ \\
				\hline
				$\{1\}$       & $\varnothing$ & $\varnothing$ & $\varnothing$ & $\varnothing$ & $\varnothing$ & $\varnothing$ & $\{0.2\}$     & $\varnothing$ & $\varnothing$ \\
				\hline
				$\varnothing$ & $\varnothing$ & $\varnothing$ & $\varnothing$ & $\varnothing$ & $\{0,1\}$     & $\{0.55\}$    & $\{1\}$       & $\varnothing$ & $\varnothing$ \\
				\hline
				$\varnothing$ & $\{0.64\}$    & $\{0.6\}$     & $\{1\}$       & $\varnothing$ & $\varnothing$ & $\varnothing$ & $\varnothing$ & $\varnothing$ & $\varnothing$ \\
				\hline
				$\{0.4\}$     & $\varnothing$ & $\{0.6\}$     & $\{1\}$       & $\{0.3\}$     & $\varnothing$ & $\varnothing$ & $\varnothing$ & $\varnothing$ & $\varnothing$ \\
				\hline
				$\varnothing$ & $\varnothing$ & $\varnothing$ & $\varnothing$ & $\varnothing$ & $\varnothing$ & $\varnothing$ & $\{0.2\}$     & $\{0.2,0.8\}$ & $\varnothing$ \\
				\hline
				$\varnothing$ & $\varnothing$ & $\{0\}$       & $\varnothing$ & $\varnothing$ & $\varnothing$ & $\varnothing$ & $\{1\}$       & $\varnothing$ & $\{0.6\}$     \\
				\hline
				$\varnothing$ & $\varnothing$ & $\varnothing$ & $\varnothing$ & $\varnothing$ & $\varnothing$ & $\varnothing$ & $\varnothing$ & $\{0.8\}$     & $\varnothing$ \\
				\hline
				$\varnothing$ & $\varnothing$ & $\{0\}$       & $\varnothing$ & $\varnothing$ & $\varnothing$ & $\varnothing$ & $\varnothing$ & $\varnothing$ & $\{0.6\}$     \\
				\hline
				$\varnothing$ & $\varnothing$ & $\varnothing$ & $\varnothing$ & $\varnothing$ & $\{0,1\}$     & $\varnothing$ & $\{1\}$       & $\varnothing$ & $\varnothing$ \\
				\hline
			\end{tabular}%
		\end{table}
	\end{step}
	
	\begin{step}\label{st-2}
		From Table \ref{t_3}, it is clear that $I_{j} \neq \varnothing$, $\forall j \in \mathscr{J}$.
	\end{step}
	
	\begin{step}\label{st-3}
		As it was described in Table \ref{t-4}, for each $i \in \mathscr{I}$ there exists at least one $j_{i} \in \mathscr{J}$ such that $S_{i j_{i}}^{\prime} \neq \varnothing$.
	\end{step}
	
	\begin{step}\label{st-4}
		Since $I_{10}=\{0.6\}$ is a singleton set (see Table \ref{t_3}) and $0.6 \in S_{7,10}^{\prime} \cap S_{9,10}^{\prime}$ (see Table \ref{t-4}), Corollary \ref{corl-8} indicates that $x_{10}=0.6$ for each feasible solution $x$ (particularly, $x_{10}^{*}=0.6$ for each optimal solution $x^{*}$ ). Also, column 10 and both rows 7 and 9 can be removed from the problem. So, by this simplification technique, the upper bound of the number of admissible functions is reduced from 5184 (Example \ref{ex-2}) to $\prod_{i \in \mathscr{I}-\{7,9\}}\left|\mathscr{J}_{i}\right|=3 \times 2 \times 3 \times 3 \times 4 \times 2 \times 1 \times 2=864$; that is $|E| \leq 864$.
		
		Moreover, since in Table \ref{t-4}, $S_{10, j}^{\prime} \subseteq S_{3 j}^{\prime}$, $\forall j \in \mathscr{J}$, it is concluded that the third equation is a redundant constraint, and it can be deleted by Corollary \ref{corl-9}. Hence, by applying this simplification technique, the upper bound of the number of admissible functions is further reduced from 864 to $\prod_{i \in \mathscr{I}-\{3,7,9\}}\left|\mathscr{J}_{i}\right|=3 \times 2 \times 3 \times 4 \times 2 \times 1 \times 2=288$; that is $|E| \leq 288$.
		
		Based on Table \ref{t-4}, it follows that the tenth row itself can also be deleted by Corollary \ref{corl-11}; because, $S_{10,6}^{\prime}=\{0,1\}$ and therefore $\left|S_{10,6}^{\prime}\right|=2$. As a result, the above-mentioned upper bound is reduced from 288 to $\prod_{i \in \mathscr{I}-\{3,7,9,10\}}\left|\mathscr{J}_{i}\right|=3 \times 2 \times 3 \times 4 \times 2 \times 1=144$, and therefore $|E| \leq 144$.
		
		Subsequently, according to Table \ref{t-4}, it turns out that $\mathscr{J}_{8}=\{9\}$ and $S_{89}^{\prime}=\{0.8\}$ are singleton sets and $0.8 \in S_{69}^{\prime}$. Thus, $x_{9}^{*}$ is assigned to $x_{9}^{*}=0.8$, and also column 9 and rows 6 and 8 can be deleted by Corollary \ref{corl-10}. Consequently, we have $\prod_{i \in \mathscr{I}-\{3,6,7,8,9,10\}}|\mathscr{J}_{i}|=3 \times 2 \times 3 \times 4=72$ and $|E| \leq 72$.
		
		So, after applying the above simplification techniques, columns $\{9,10\}$ and rows $\{3,6,7,8,9,10\}$ are deleted and we obtain $x_{9}^{*}=0.8$ and $x_{10}^{*}=0.6$. Therefore, the reduced matrices $A^{+}$ and $A^{-}$, and the right-hand-side vector $b$ become
		
		$$A^{+}=\left[\begin{array}{llllllll}
			0.25 & 0.32 & 0.41 & 0.19 & 0.70 & 0.13 & 0.44 & 0.37 \\
			0.80 & 0.73 & 0.64 & 0.79 & 0.80 & 0.22 & 0.80 & 0.56 \\
			0.10 & 0.23 & 0.25 & 0.15 & 0.12 & 0.05 & 0.02 & 0.01 \\
			0.45 & 0.35 & 0.70 & 0.50 & 0.41 & 0.27 & 0.39 & 0.48
		\end{array}\right]$$
		
		$$A^{-}=\left[\begin{array}{llllllll}
			0.70 & 0.70 & 0.32 & 0.44 & 0.00 & 0.16 & 0.20 & 0.50 \\
			0.70 & 0.65 & 0.14 & 0.12 & 0.80 & 0.76 & 0.00 & 1.00 \\
			0.14 & 0.10 & 0.04 & 0.00 & 0.10 & 0.00 & 0.14 & 0.02 \\
			0.70 & 0.04 & 0.27 & 0.36 & 0.60 & 0.40 & 0.48 & 0.50
		\end{array}\right]$$
		
		$$
		\begin{aligned}
			& b^{T}=[0.50,0.80,0.15,0.50]
		\end{aligned}
		$$
		The current matrices $A^{+}$ and $A^{-}$ are equivalent to four rows (rows 1,2,4 and 5 in the main problem) and eight columns 1-8. Furthermore, Table \ref{t-4} is updated as follows:
		
		\begin{table}[!ht]
			\centering
			\caption{Updated Table 4 after applying Corollaries 8 - 11.}
			\label{t-5}
			\begin{tabular}{|c|c|c|c|c|c|c|c|}
				\hline
				$\{0.4\}$     & $\{0.4\}$     & $\varnothing$ & $\varnothing$ & $\{0.6\}$     & $\varnothing$ & $\varnothing$ & $\varnothing$ \\
				\hline
				$\{1\}$       & $\varnothing$ & $\varnothing$ & $\varnothing$ & $\{1\}$       & $\varnothing$ & $\varnothing$ & $\{0.2\}$     \\
				\hline
				$\varnothing$ & $\{0.64\}$    & $\{0.6\}$     & $\{1\}$       & $\varnothing$ & $\varnothing$ & $\varnothing$ & $\varnothing$ \\
				\hline
				$\{0.4\}$     & $\varnothing$ & $\{0.6\}$     & $\{1\}$       & $\{0.3\}$     & $\varnothing$ & $\varnothing$ & $\varnothing$ \\
				\hline
			\end{tabular}%
		\end{table}
		
		By considering Table \ref{t-5}, it follows that $\mathscr{I}_{8}=\{0.2\}$ and $S_{28}^{\prime}=\{0.2\}$. So, we have $L_{8}=0.2$ (see Table \ref{t_3}) and therefore $L_{8} \in S_{i 8}^{\prime}$, $\forall i \in \mathscr{I}_{8}$. Hence, Corollary \ref{corl-12} implies that $x(e)_{8}=L_{8}=0.2$, $\forall e \in E$. Consequently, by Theorem \ref{thm-4}, it is concluded that $x_{8}^{*}=0.2$ for each optimal solution $x^{*}$. Also, column 8 and row 2 can be removed from the problem. So, this simplification technique leads us to $\prod_{i \in \mathscr{I}-\{2,3,6,7,8,9,10\}}\left|\mathscr{J}_{i}\right|=3 \times 3 \times 4=36$ and $|E| \leq 36$. In addition, since $S_{i 6}^{\prime}=S_{i 7}^{\prime}=\varnothing$ for each $i \in\{1,2,3,4\}$ (see Table \ref{t-5}), from Definition \ref{def-5} we have $6 \notin \mathscr{J}_{i}$ and $7 \notin \mathscr{J}_{i}$, $\forall i \in\{1,2,3,4\}$. Thus, from Definition \ref{def-6}, $\mathscr{I}_{6}(e)=\mathscr{I}_{7}(e)=\varnothing$, $\forall e \in E$, which together with (\ref{eq-12}) imply $x(e)_{6}=L_{6}$ and $x(e)_{7}=L_{7}$, $\forall e \in E$ (where $L_{6}=L_{7}=0$ from Table \ref{t_3}). So, by using Theorem \ref{thm-4}, we can set $x_{6}^{*}=x_{7}^{*}=0$ and delete columns 6 and 7 . After these reductions, matrices $A^{+}$and $A^{-}$, vector $b$ and Table \ref{t-5} are updated as follows:
		
		$$
		\begin{gathered}
			A^{+}=\left[\begin{array}{lllll}
				0.25 & 0.32 & 0.41 & 0.19 & 0.70 \\
				0.10 & 0.23 & 0.25 & 0.15 & 0.12 \\
				0.45 & 0.35 & 0.70 & 0.50 & 0.41
			\end{array}\right] \quad A^{-}=\left[\begin{array}{lllll}
				0.70 & 0.70 & 0.32 & 0.44 & 0.00 \\
				0.14 & 0.10 & 0.04 & 0.00 & 0.10 \\
				0.70 & 0.04 & 0.27 & 0.36 & 0.60
			\end{array}\right] \\
			b^{T}=[0.50,0.15,0.50]
		\end{gathered}
		$$
		
		\begin{table}[!ht]
			\centering
			\caption{Updated Table \ref{t-5} after applying Corollary \ref{corl-12}.}
			\label{t-6}
			\begin{tabular}{|c|c|c|c|c|}
				\hline
				$\{0.4\}$     & $\{0.4\}$     & $\varnothing$ & $\varnothing$ & $\{0.6\}$     \\
				\hline
				$\varnothing$ & $\{0.64\}$    & $\{0.6\}$     & $\{1\}$       & $\varnothing$ \\
				\hline
				$\{0.4\}$     & $\varnothing$ & $\{0.6\}$     & $\{1\}$       & $\{0.3\}$     \\
				\hline
			\end{tabular}%
		\end{table}
		
		The situation of Table \ref{t-6} and the reduced matrices $A^{+}$ and $A^{-}$are the same as having three rows (rows 1,4 and 5 in the main problem) and five columns 1-5. As it is shown in Table \ref{t-6}, $\mathscr{I}_{1}=\mathscr{I}_{5}=\{1,3\}$ and $\bigcap_{i \in \mathscr{I}_{1}} S_{i 1}^{\prime}=S_{11}^{\prime} \cap S_{31}^{\prime}=\{0.4\}$. On the other hand, $L_{1}=0.4$ (see Table \ref{t_3}), and therefore Corollary \ref{corl-13} (Part (a)) implies that $x_{5}^{*}=L_{5}=0.3$. Also, column 5 can be removed from the problem. Similarly, since $\mathscr{I}_{3}=\mathscr{I}_{4}=\{2,3\}$, $\bigcap_{i \in \mathscr{I}_{3}} S_{i 3}^{\prime}=\{0.6\}=\left\{U_{3}\right\}$, $\bigcap_{i \in \mathscr{I}_{4}} S_{i 4}^{\prime}=\{1\}=\left\{U_{4}\right\}$ and $c_{3}\left(U_{3}-L_{3}\right)=0.93(0.6-0)=0.558<3.28=3.28(1-0)=c_{4}\left(U_{4}-L_{4}\right)$, from Part (b) of Corollary \ref{corl-13} , we can set $x_{4}^{*}=L_{4}=0$ and remove column 4 . By applying this simplification technique, for the remaining rows (rows 1, 4 and 5 in the main problem), sets $\mathscr{J}_{1}=\{1,2,5\}$, $\mathscr{J}_{4}=\{2,3,4\}$ and $\mathscr{J}_{5}=\{1,3,4,5\}$ are updated as $\mathscr{J}_{1}=\{1,2\}$, $\mathscr{J}_{4}=\{2,3\}$ and $\mathscr{J}_{5}=\{1,3\}$, respectively. So, we have $\prod_{i \in \mathscr{I}-\{2,3,6,7,8,9,10\}}\left|\mathscr{J}_{i}\right|=2 \times 2 \times 2=8$ and reduced matrices $A^{+}$ and $A^{-}$, vector $b$ and Table \ref{t-6} are updated as follows:
		$$
		\begin{gathered}
			A^{+}=\left[\begin{array}{lll}
				0.25 & 0.32 & 0.41 \\
				0.10 & 0.23 & 0.25 \\
				0.45 & 0.35 & 0.70
			\end{array}\right] \quad A^{-}=\left[\begin{array}{lll}
				0.70 & 0.70 & 0.32 \\
				0.14 & 0.10 & 0.04 \\
				0.70 & 0.04 & 0.27
			\end{array}\right] \\
			b^{T}=[0.50,0.15,0.50]
		\end{gathered}
		$$
		\begin{table}[!ht]
			\centering
			\caption{Updated Table 6 after applying Corollary \ref{corl-13}.}
			\label{t-7}
			\begin{tabular}{|c|c|c|}
				\hline
				$\{0.4\}$     & $\{0.4\}$     & $\varnothing$ \\
				\hline
				$\varnothing$ & $\{0.64\}$    & $\{0.6\}$     \\
				\hline
				$\{0.4\}$     & $\varnothing$ & $\{0.6\}$     \\
				\hline
			\end{tabular}%
		\end{table}
		
		Table \ref{t-6} and the reduced matrices $A^{+}$ and $A^{-}$ are equivalent to three rows (rows 1, 4 and 5 in the main problem) and three columns 1-3 associated with $x_{1}$, $x_{2}$ and $x_{3}$, respectively.
	\end{step}
	
	\begin{step}\label{st-5}
		Figure \ref{fig_1} contains the solution tree generated by the branch and bound approach, which is based on the concept of modified functions. This method begins by choosing $e(1)$ from $\bar{\mathscr{J}}_{1}(e)=\mathscr{J}_{1}=\{1,2\}$. Hence, each $e \in \bar{E}$ must satisfy either $e(1)=1$ or $e(1)=2$. This yields two branches denoted by nodes 1 and 2 in Figure \ref{fig_1} . For node 1, by considering $e_{1} \in \bar{E}$ and $e_{1}(1)=1$, we have $e_{1}=[1]$, $\mathscr{I}_{1}\left(e_{1}\right)=\{1\}$ and $\mathscr{I}_{2}\left(e_{1}\right)=\mathscr{I}_{3}\left(e_{1}\right)=\varnothing$. So, from (\ref{eq-12}) we obtain $x\left(e_{1}\right)_{1}=\min \left\{\bigcap_{i \in \mathscr{I}_{1}\left(e_{1}\right)} S_{i 1}^{\prime}\right\}=\min S_{11}^{\prime}=0.4$, $x\left(e_{1}\right)_{2}=L_{2}=0.4$ and $x\left(e_{1}\right)_{3}=L_{3}=0$ with $z_{1}=c_{1} x\left(e_{1}\right)_{1}+c_{2} x\left(e_{1}\right)_{2}+c_{3} x\left(e_{1}\right)_{3}=1(0.4)+0.35(0.4)+0.93(0)=0.54$. Similarly, by selecting $e(1)=2$ and considering another admissible function $e_{2} \in \bar{E}$ such that $e_{2}(1)=2$, it follows that $e_{2}=[2]$, $x\left(e_{2}\right)=[0.4,0.4,0]$ and $z_{2}=0.54$. At this time, we have no reason to exclude nodes 1 and 2 from consideration; therefore, both nodes should be further investigated. By selecting node 1 (where $e_{1}(1)=1$ ), for the second row (row 4 in the main problem), we have $\mathscr{J}_{4}=\{2,3\}$. Since $\mathscr{I}_{2}\left(e_{1}, 2\right)=\mathscr{I}_{3}\left(e_{1}, 2\right)=\varnothing$, then $\mathscr{J}_{2}\left(e_{1}\right)=\{2,3\}$ (Definition \ref{def-6}). Also, since $\mathscr{J}_{2}\left(e_{1}\right) \cap\left\{e_{1}(1)\right\}=\varnothing$, we have $\bar{\mathscr{J}}_{2}\left(e_{1}\right)=\mathscr{J}_{2}\left(e_{1}\right)=\{2,3\}$ (Definition \ref{def-9}). This implies that for each $e \in \bar{E}$ such that $e(1)=1$ , we must have either $e(2)=2$ or $e(2)=3$. So, branching on node 1 yields nodes 3 and 4 in Figure \ref{fig_1}, which lead us to two admissible functions $e_{1}$ and $e_{3}$, respectively. At node 3, $e_{1}(2)=2$ and therefore we obtain $e_{1}=[1,2]$, $x\left(e_{1}\right)=[0.4,0.64,0]$ (resulted by (\ref{eq-12})) and $z_{1}=c_{1} x\left(e_{1}\right)_{1}+c_{2} x\left(e_{1}\right)_{2}+c_{3} x\left(e_{1}\right)_{3}=1(0.4)+0.35(0.64)+0.93(0)=0.624$. On the other hand, if we move along $e(2)=3$ to node 4 , similar reasoning shows that $e_{3}=[1,3]$, $x\left(e_{3}\right)=[0.4,0.4,0.6]$ and $z_{3}=1.098$ for node 4. Again, there is no reason to exclude nodes 3 and 4 from consideration, so we need to branch on one node. By using the jump-tracking technique, a branch is made at a node that has a lower bound on $z$. Here, node 3 has been chosen.
		
		At node 3, $e_{1}(1)=1$ and $e_{1}(2)=2$. For the third row (row 5 in the main problem), we have $\mathscr{J}_{5}=\{1,3\}$. Since $\mathscr{I}_{1}\left(e_{1}, 3\right)=\{1\}$, $S_{31}^{\prime} \cap\left(\bigcap_{k \in \mathscr{I}_{1}\left(e_{1}, 3\right)} S_{k 1}^{\prime}\right)=S_{31}^{\prime} \cap S_{11}^{\prime}=\{0.4\}$ and $\mathscr{I}_{3}\left(e_{1}, 3\right)=\varnothing$, then $\mathscr{J}_{3}\left(e_{1}\right)=\{1,3\} \quad$ (Definition \ref{def-6}). However, $\mathscr{J}_{3}\left(e_{1}\right) \cap\left\{e_{1}(1), e_{1}(2)\right\}=\mathscr{J}_{3}\left(e_{1}\right) \cap\{1,2\}=\{1\}$, which means $\bar{\mathscr{J}}_{3}\left(e_{1}\right)=\{1\}$ (Definition \ref{def-9}). This implies that $e_{1}(3)=1$ is the only selection for constructing a modified function. So, any choice associated with node 3 must satisfy $e_{1}(3)=1$, which yields one branch from node 3, namely, node 5 in Figure \ref{fig_1}. By the same reasoning as before, for this node, we have $e_{1}=[1,2,1]$, $x\left(e_{1}\right)=[0.4,0.64,0]$ and $z_{1}=0.624$. Node 5 corresponds to the complete branch $e_{1}=[1,2,1]$ and $x\left(e_{1}\right)=[0.4,0.64,0]$ may be viewed as a candidate solution with $z_{1}=0.624$.
		
		Because the objective value for node 4 ($z=1.098$) cannot be lower than 0.624, this node can be eliminated from further consideration. Hence, we now branch on node 2. At node 2, $e_{2}(1)=2$. So, at the second row (row 4 in the main problem), we have $\mathscr{J}_{4}=\{2,3\}$. Since $\mathscr{I}_{2}\left(e_{2}, 2\right)=\{1\}$, $S_{22}^{\prime} \cap\left(\bigcap_{k \in \mathscr{I}_{2}\left(e_{2}, 2\right)} S_{k 2}^{\prime}\right)=S_{22}^{\prime} \cap S_{12}^{\prime}=\{0.6\} \cap\{0.4\}=\varnothing$ and $\mathscr{I}_{3}\left(e_{2}, 2\right)=\varnothing$, then $\mathscr{J}_{2}\left(e_{2}\right)=\{3\}$ (Definition \ref{def-6}). Also, $\mathscr{J}_{2}\left(e_{2}\right) \cap\left\{e_{2}(1)\right\}=\mathscr{J}_{2}\left(e_{2}\right) \cap\{2\}=\varnothing$, which means $\bar{\mathscr{J}}_{2}\left(e_{2}\right)=\mathscr{J}_{2}\left(e_{2}\right)=\{3\}$ (Definition \ref{def-9}). This implies that $e_{2}(2)=3$ is the only selection for constructing a modified function. So, any modified function associated with node 2 must satisfy $e_{2}(2)=3$, that yields one branch from node 2 . Correspondingly, we create node 6 in Figure \ref{fig_1}, where $e_{2}=[2,3]$, $x\left(e_{2}\right)=[0.4,0.4,0.6]$ and $z_{2}=1.098$. Thus, since $z_{2}>z_{1}$, this node is eliminated from consideration. As a result, the branch and bound method finds the optimal variables $x_{1}^{*}=0.4$, $x_{2}^{*}=0.64$ and $x_{3}^{*}=0$ (associated with the solution $x\left(e_{1}\right)=[0.4,0.64,0]$ ) by checking only one complete modified function $e_{1}=[1,2,1]$.
		\begin{figure}[!ht]
			\begin{center}
				\includegraphics[scale=0.6]{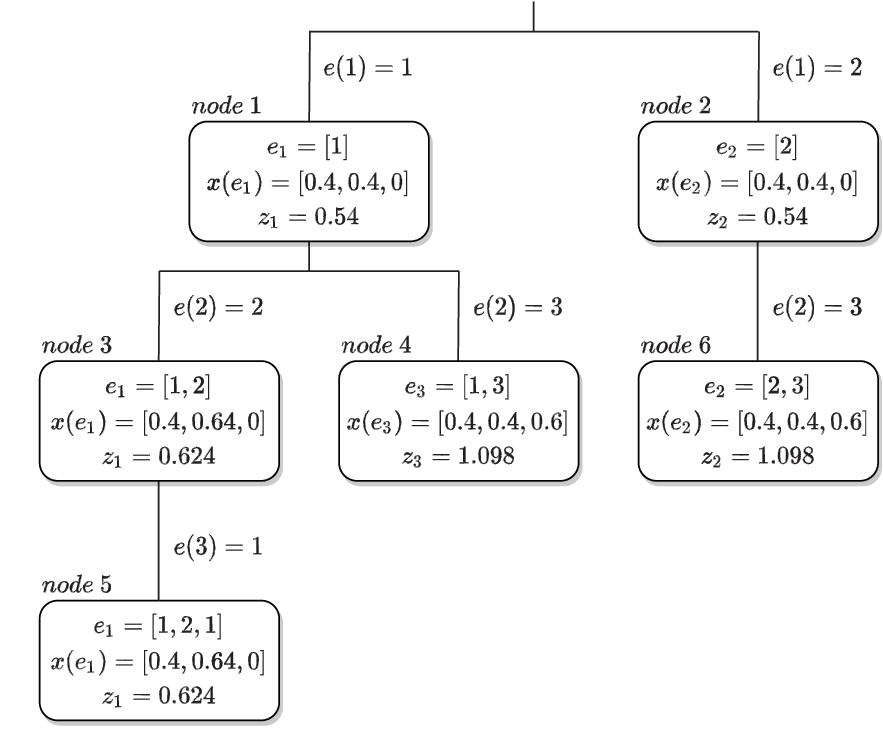}
				\captionsetup{justification=centering}
				\caption{Branch and bound method based on the modified functions}
				\label{fig_1}
			\end{center}
		\end{figure}
	\end{step}
	\begin{step}\label{st-6}
		According to Steps \ref{st-4} and \ref{st-5}, an optimal solution of the problem is obtained as $$x^{*}=(0.4,0.64,0,0,0.3,0,0,0.2,0.8,0.6)$$ with $cx^{*}=10.717$.
	\end{step}
\end{example}


\section*{Conclusion}\label{sec_con}
This paper proposes an algorithm to find a global optimal solution to linear problems involving a wide class of bipolar fuzzy relation equations defined with continuous Archimedean t-norms (including both continuous strict and nilpotent t-norms). There were several basic properties of continuous Archimedean bipolar equations explored in this study, and a finite number of compact sets characterized the feasible solution set. Additionally, two necessary feasibility conditions were presented as part of the feasibility determination. As well, seven simplification techniques were introduced to facilitate the process of solution finding; four techniques for deleting redundant constraints and/or eliminating redundant columns by assigning a fixed value to their corresponding variables, and three additional techniques for reducing the domain search to a subset of the feasible region, which contains the optimal solutions. Then, a branch-and-bound-based algorithm was presented that takes advantage of the concept of modified functions to accelerate the solution process by neglecting some candidate solutions that are not optimal and removing some branches that do not lead to an optimal solution. It is worth pointing out that it is a future study to consider the assumption of left-hand continuity throughout the paper, and the comparison result is well worth assessing. Despite this, since the complexity of solving the problem increases exponentially with the problem size, it is always a challenge to find other simplification rules or methods that can reach the optimal solution.



\section*{Appendix A}\label{sec_app}
\setcounter{table}{0} \renewcommand{\thetable}{A\arabic{table}}

\SetTblrInner{rowsep=3pt,colsep=3pt}
\begin{table}[!ht]
	\centering
	\caption{Summary of some commonly used continuous Archimedean t-norms.}
	\label{tbl_a1}
		\scalebox{0.8}{
			\begin{tblr}{c|c|c}
				\textbf{t-norm}  & \textbf{Function}                                                                                                                     & \textbf{Parameter} \\ \hline \hline
				Product          & $T_P(x, y)=xy$                                                                                                                        & -                  \\ \hline
				Einstein Product & $T_{EP}(x, y)=\frac{x y}{2-(x+y-x y)}$                                                                                                & -                  \\ \hline
				Lukasiewicz      & $T_{L}(x, y)=max\{0,x+y-1\}$                                                                                                          & -                  \\ \hline
				Frank            & $T_F^S(x, y)=\log _s\left(1+\frac{\left(s^x-1\right)\left(s^y-1\right)}{s-1}\right)$                                                  & $s>0, s\neq 1$     \\ \hline
				Yager            & $T_Y^P(x, y)=max\{0,1-[(1-x)^p+(1-y)^p]^{1/p}\}$                                                                                      & $P>0$              \\ \hline
				Hamacher         & $
				T_H^\alpha(x, y)=
				\begin{cases}
					0                                      & , \quad \alpha=x=y=0        \\
					\frac{x y}{\alpha+(1-\alpha)(x+y-x y)} & , \quad \text { otherwise }
				\end{cases}
				$                & $\alpha\geq 0$                                                                                                                                             \\ \hline
				Dombi            & $T_D^\lambda(x, y)=
				\begin{cases}
					0                                                                                                            & , \quad x=0 \text { or } y=0 \\
					\frac{1}{1+\left[\left(\frac{1-x}{x}\right)^\lambda+\left(\frac{1-y}{y}\right)^\lambda\right]^{1 / \lambda}} & , \quad \text {otherwise}
				\end{cases}
				$                & $\lambda>0$                                                                                                                                                \\ \hline
				Schweizer-Sklar  & $T_{S S}^p(x, y)=\left(\max \left\{0, x^p+y^p-1\right\}\right)^{1 / p}$                                                               & $p \neq 0$              \\ \hline
				Sugeno-Weber     & $T_{SW}^{\lambda}(x,y)=max\{0,(x+y-1+\lambda xy)/(1+ \lambda)\}$                                                                      & $\lambda>-1$       \\ \hline
				Aczel-Alsina     & $T_{A A}^\lambda(x, y)=\exp \left(-\left[(-\operatorname{Ln}(x))^\lambda+(-\operatorname{Ln}(y))^\lambda\right]^{1 / \lambda}\right)$ & $\lambda>0$
			\end{tblr}%
		}
	\end{table}
	
	\begin{table}[!ht]
		\centering
		\caption{Additive generator $f_{\varphi}$ and pseudoinverse $f_{\varphi}^{(-1)}$ discussed in Lemma \ref{lm-2} for the t-norms stated in Table \ref{tbl_a1}.}
		\label{tbl_a2}
		\scalebox{0.8}{
			\begin{tabular}{c|c|c|c}
				\textbf{t-norm}    & \textbf{$f_{\varphi}, f_{\varphi}^{(-1)}$}                                                                                                 & \textbf{t-norm}    & \textbf{$f_{\varphi}, f_{\varphi}^{(-1)}$}                        \\
				\hline \hline
				Product            & \begin{tabular}{l}
					$f_{P}(x)=-\operatorname{Ln}(x)$ \\
					$f_{P}^{(-1)}(x)=e^{-x}$         \\
				\end{tabular}                                                                           & Lukasiewicz        & \begin{tabular}{l}
					$f_{L}(x)=1-x$                    \\
					$f_{L}^{(-1)}(x)=\max \{1-x, 0\}$ \\
				\end{tabular}                                                          \\
				\hline
				\begin{tabular}{l}
					Einstein \\
					Product  \\
				\end{tabular} & \begin{tabular}{l}
					$f_{E P}(x)=\operatorname{Ln}((2-x) / x)$   \\
					$f_{E P}^{(-1)}(x)=2 /\left(1+e^{x}\right)$ \\
				\end{tabular}                                                                     & Frank              & \begin{tabular}{l}
					$f_{F}(x)=\log _{s}\left((s-1) /\left(s^{x}-1\right)\right)$         \\
					$f_{F}^{(-1)}(x)=\log _{s}\left(\left((s-1) / s^{x}\right)+1\right)$ \\
				\end{tabular} \\
				\hline
				Yager              & \begin{tabular}{l}
					$f_{Y}(x)=(1-x)^{p}$                        \\
					$f_{Y}^{(-1)}(x)=\max \{1-\sqrt[p]{x}, 0\}$ \\
				\end{tabular}                              & \begin{tabular}{l}
					Sugeno- \\
					Weber   \\
				\end{tabular} & \begin{tabular}{l}
					$f_{S W}(x)=1-\log _{1+\lambda}(1+\lambda x)$                                         \\
					$f_{S W}^{(-1)}(x)=\max \left\{\left((1+\lambda)^{1-x}-1\right) / \lambda, 0\right\}$ \\
				\end{tabular}                              \\
				\hline
				Dombi              & \begin{tabular}{l}
					$f_{D}(x)=((1-x) / x)^{\lambda}$           \\
					$f_{D}^{(-1)}(x)=1 /(1+\sqrt[\lambda]{x})$ \\
				\end{tabular}                                     & Aczel-Alsina       & \begin{tabular}{l}
					$f_{A A}(x)=(-\operatorname{Ln}(x))^{\lambda}$ \\
					$f_{A A}^{(-1)}(x)=e^{-\sqrt[\lambda]{x}}$     \\
				\end{tabular}                                                   \\
				\hline
				\begin{tabular}{l}
					Schweizer- \\
					Sklar      \\
				\end{tabular} & \begin{tabular}{l}
					$f_{S S}(x)=\left(1-x^{p}\right) / p$                                                                       \\
					$f_{S S}^{(-1)}(x)= \begin{cases}\sqrt[p]{1-p x} & , p<0 \\ \max \{\sqrt[p]{1-p x}, 0\} & , p>0\end{cases}$ \\
				\end{tabular} & Hamacher           & \begin{tabular}{l}
					$f_{H}(x)= \begin{cases}(1-x) / x & , \alpha=0 \\ \operatorname{Ln}((\alpha+(1-\alpha) x) / x) & , \alpha>0\end{cases}$ \\
					$f_{H}^{(-1)}(x)= \begin{cases}1 /(1+x) & , \alpha=0 \\ \alpha /\left(\alpha-1+e^{x}\right) & , \alpha>0\end{cases}$    \\
				\end{tabular}                                                                                                  \\
		\end{tabular}}
	\end{table}
	
	\begin{table}[!ht]
		\centering
		\caption{Value $u$ discussed in Theorems \ref{thm-1} and \ref{thm-2} for the t-norms stated in Table \ref{tbl_a1}.}
		\label{tbl_a3}
		\resizebox{\columnwidth}{!}{%
			\scalebox{0.8}{
				\begin{tabular}{c|c|c|c}
					\textbf{t-norm}  & \textbf{$u$}                                   & \textbf{t-norm} & \textbf{$u$}                                                                             \\ \hline \hline
					Product          & $b / a$                                      & Lukasiewicz     & $1+b-a$                                                                                \\ \hline
					Einstein Product & $((2-a) b) /(a+b-a b)$                       & Frank           & $\log _{s}\left(1+\left[\left(s^{b}-1\right)(s-1) /\left(s^{a}-1\right)\right]\right)$ \\ \hline
					Yager            & $1-\left((1-b)^{p}-(1-a)^{p}\right)^{1 / p}$ & Sugeno-Weber    & $((1+\lambda) b+1-a) /(1+\lambda a)$                                                   \\ \hline
					Dombi &
					$\left(1+\left[((1-b) / b)^{\lambda}-((1-a) / a)^{\lambda}\right]^{1 / \lambda}\right)^{-1}$ &
					Aczel-Alsina &
					$\exp \left(-\left[(-\operatorname{Ln}(b))^{\lambda}-(-\operatorname{Ln}(a))^{\lambda}\right]^{1 / \lambda}\right)$ \\ \hline
					Schweizer-Sklar  & $\left(1+b^{p}-a^{p}\right)^{1 / p}$         & Hamacher        & $([\alpha+(1-\alpha) a] b) /(a-(1-\alpha)(1-a) b)$                                     \\ 
			\end{tabular}}%
		}
	\end{table}
	
	\newpage
	
	
	\section*{Appendix B}\label{sec_appB}
	\begin{proof}[The proof of Corollary \ref{corl-6}]
		Suppose that $e \in E$ and $\mathscr{I}_{j_{0}}(e) \neq \varnothing$ for some $j_{0} \in \mathscr{J}$. Also, without loss of generality, let $i_{0}=\max \mathscr{I}_{j_{0}}(e)$. So, according to Definition \ref{def-6} , we have $\mathscr{I}_{j_{0}}(e)=\mathscr{I}_{j_{0}}\left(e, i_{0}\right) \cup\left\{i_{0}\right\}$ and $S_{i_{0} j_{0}}^{\prime} \cap\left(\bigcap_{k \in \mathscr{I}_{J_{0}}\left(e, i_{0}\right)} S_{kj_{0}}^{\prime}\right) \neq \varnothing$. Thus, $\bigcap_{i \in \mathscr{I}_{j_{0}}(e)} S_{ij_{0}}^{\prime} \neq \varnothing$. To prove the converse statement, since $e(i) \in \mathscr{J}_{i}$ ($\forall i \in \mathscr{I}$), it is sufficient to show that if $e\left(i_{0}\right)=j_{0}$ and $\mathscr{I}_{j_{0}}\left(e, i_{0}\right) \neq \varnothing$, then $S_{i_{0} j_{0}}^{\prime} \cap\left(\bigcap_{k \in \mathscr{I}_{j_{0}}\left(e, i_{0}\right)} S_{kj_{0}}^{\prime}\right) \neq \varnothing$. However, since $\mathscr{I}_{j_{0}}\left(e, i_{0}\right) \subseteq \mathscr{I}_{j_{0}}(e)$, the statement $\mathscr{I}_{j_{0}}\left(e, i_{0}\right) \neq \varnothing$ implies $\mathscr{I}_{j_{0}}(e) \neq \varnothing$. Hence, from the assumption, we have $\bigcap_{i \in \mathscr{I}_{j_{0}}(e)} S_{ij_{0}}^{\prime} \neq \varnothing$. Therefore, since $i_{0} \in \mathscr{I}_{j_{0}}(e)$ (because $\left.e\left(i_{0}\right)=j_{0}\right)$ and $\mathscr{I}_{j_{0}}\left(e, i_{0}\right) \subseteq \mathscr{I}_{j_{0}}(e)$, it is concluded that $S_{i_{0} j_{0}}^{\prime} \cap\left(\bigcap_{k \in \mathscr{I}_{j_0}\left(e, i_{0}\right)} S_{kj_{0}}^{\prime}\right) \neq \varnothing$.
	\end{proof}
	
	\begin{proof}[The proof of Theorem \ref{thm-3}]
		Let $x \in \bigcup_{e \in E} S(e)$. So, $x \in S\left(e_{0}\right)$ for some $e_{0} \in E$. Hence, according to (\ref{eq-11}), for each $j \in \mathscr{J}$ we have either $x_{j} \in \bigcap_{i \in \mathscr{I}_{j}\left(e_{0}\right)} S_{i j}^{\prime}$ (if $\mathscr{I}_{j}\left(e_{0}\right) \neq \varnothing$ ) or $x_{j} \in I_{j}$ (if $\mathscr{I}_{j}\left(e_{0}\right)=\varnothing$ ). But, since $S_{i j}^{\prime}=S_{i j} \cap I_{j}$ (Definition \ref{def-4}), from $x_{j} \in \bigcap_{i \in \mathscr{I}_{j}\left(e_{0}\right)} S_{i j}^{\prime}$ it is obtained again $x_{j} \in I_{j}$. Consequently, $x_{j} \in I_{j}$, $\forall j \in \mathscr{J}$. On the other hand, from Definition \ref{def-6} we have $e_{0}(i)=j_{i} \in \mathscr{J}_{i}\left(e_{0}\right) \subseteq \mathscr{J}_{i}$ $(\forall i \in \mathscr{I})$, which implies $\mathscr{I}_{j_{i}}\left(e_{0}\right) \neq \varnothing$, and therefore from (\ref{eq-11}) we have $x_{j_{i}} \in \bigcap_{k \in \mathscr{I}_{j_{i}}\left(e_{0}\right)} S_{kj_{i}}^{\prime} \subseteq S_{i j_{i}}^{\prime}$. Hence, $x_{j_{i}} \in S_{i j_{i}}^{\prime}$, $\forall i \in \mathscr{I}$. Now, Lemma \ref{lm-5} requires that $x \in S\left(A^{+}, A^{-}, b\right)$. Conversely, let $x \in S\left(A^{+}, A^{-}, b\right),~ \mathscr{I}_{j}(x)=\left\{i \in \mathscr{I}: x_{j} \in S_{i j}^{\prime}\right\}$ and $\mathscr{J}_{i}(x)=\left\{j \in \mathscr{J}: x_{j} \in S_{i j}^{\prime}\right\}$. So, for each $i \in \mathscr{I}_{j}(x)$ and each $j \in \mathscr{J}_{i}(x)$, we have $S_{i j}^{\prime} \neq \varnothing$ that means $j \in \mathscr{J}_{i}$ (Definition \ref{def-5}). Also, Lemma \ref{lm-5} implies that $\mathscr{J}_{i}(x) \neq \varnothing$, $\forall i \in \mathscr{I}$. Without loss of generality, let $j_{i}=\min \mathscr{J}_{i}(x)$ and $e_{0}(i)=j_{i}$, $\forall i \in \mathscr{I}$. Therefore, $e_{0}$ is a function on $\mathscr{I}$ such that
		\begin{equation}\label{eq-13}
			e_{0}(i)=j_{i} \in \mathscr{J}_{i}, \forall i \in \mathscr{I}
		\end{equation}
		Moreover, we have
		\begin{equation}\label{eq-14}
			\mathscr{I}_{j}\left(e_{0}\right) \neq \varnothing, \forall j \in\left\{j_{1}, \ldots, j_{m}\right\}~ \text{and} ~\mathscr{I}_{j}\left(e_{0}\right)=\varnothing, \forall j \in \mathscr{J}-\left\{j_{1}, \ldots, j_{m}\right\}
		\end{equation}
		In addition, if $j_{p} \in\left\{j_{1}, \ldots, j_{m}\right\}$ and $k \in \mathscr{I}_{j_{p}}\left(e_{0}\right)$ (i.e., $\left.e_{0}(k)=j_{p}\right)$, then by our definition it follows that $j_{p}=\min \mathscr{J}_{k}(x)$ which means $x_{j_{p}} \in S_{k j_{p}}^{\prime}$. Hence, it is concluded that
		\begin{equation}\label{eq-15}
			x_{j} \in \bigcap_{i \in \mathscr{I}_{j}\left(e_{0}\right)} S_{i j}^{\prime}, \forall j \in\left\{j_{1}, \ldots, j_{m}\right\}
		\end{equation}
		Consequently, by Corollary \ref{corl-6} and (\ref{eq-12}) - (\ref{eq-14}) we have $e_{0} \in E$. Also, since $x \in S\left(A^{+}, A^{-}, b\right)$, then $x_{j} \in I_{j}$, $\forall j \in \mathscr{J}$ (Lemma \ref{lm-5}). Particularly, $x_{j} \in I_{j}$, if $\mathscr{I}_{j}\left(e_{0}\right)=\varnothing$. This fact together with (\ref{eq-11}), (\ref{eq-13}) and (\ref{eq-14}) imply $x \in S\left(e_{0}\right)$.
	\end{proof}	
	
	\begin{proof}[The proof of Lemma \ref{lm-6}]
		From Theorem \ref{thm-3}, $S_{\left\{i_{0}\right\}}\left(A^{+}, A^{-}, b\right)=\bigcup_{e^{\prime} \in E^{\prime}} S\left(e^{\prime}\right)$, where $E^{\prime}$ is the set of all the restrictions of admissible functions $e \in E$ to $\mathscr{I}-\left\{i_{0}\right\}$. So, for each $x \in S_{\left\{i_{0}\right\}}\left(A^{+}, A^{-}, b\right)$, there exists at least one $e_{0}^{\prime} \in E^{\prime}$ such that $x \in S\left(e_{0}^{\prime}\right)$. Now, from (\ref{eq-11}), for each $j \in \mathscr{J}$ we have either $x_{j} \in \bigcap_{i \in \mathscr{I}_{j}\left(e_{0}^{\prime}\right)} S_{i j}^{\prime}$ or $x_{j} \in I_{j}$. But, in the former case, since $S_{i j}^{\prime}=S_{i j} \cap I_{j}$ (Definition \ref{def-4}), we have again $x_{j} \in \bigcap_{i \in \mathscr{I}_{j}\left(e_{0}^{\prime}\right)} S_{i j}^{\prime} \subseteq I_{j}$. Consequently, $x_{j} \in I_{j}$ (and therefore $x_{j} \in I_{i_{0} j}$ from Definition \ref{def-4}), $\forall j \in \mathscr{J}$.
	\end{proof}
	
	\begin{proof}[The proof of Lemma \ref{lm-7}]
		To prove the lemma, it is sufficient to show that $S_{\left\{i_{0}\right\}}\left(A^{+}, A^{-}, b\right) \subseteq S\left(A^{+}, A^{-}, b\right)$. In other words, we shall show that if $x \in S_{\left\{i_{0}\right\}}\left(A^{+}, A^{-}, b\right)$, then $x$ satisfies the $i_{0}$'th equation. For this purpose, we first note from Lemma \ref{lm-6} that $x_{j} \in I_{j}$, $\forall j \in \mathscr{J}$ $\left(*1\right)$ and particularly $x_{j} \in I_{i_{0} j}$, $\forall j \in \mathscr{J}$ $\left(*2\right)$. Also, since $S\left(A^{+}, A^{-}, b\right) \neq \varnothing$, there exists at least one $j_{0} \in \mathscr{J}$ such that $S_{i_{0} j_{0}}^{\prime} \neq \varnothing$ (Lemma \ref{lm-4} (b)). So, from $b_{i_{0}}=0$ and Remark \ref{rmk-2}, we have $S_{i_{0} j_{0}}^{\prime}=\left[1-u_{i_{0} j_{0}}^{-}, u_{i_{0} j_{0}}^{+}\right]$ which together with Corollary \ref{corl-4} imply $S_{i_{0} j_{0}}^{\prime}=\left[L_{j_{0}}, U_{j_{0}}\right]=I_{j_{0}} (*3)$. Hence, from $\left(*1\right)$ and $\left(*3\right)$, it follows that $x_{j_{0}} \in S_{i_{0} j_{0}}^{\prime}\left(*4\right)$. Now, the result follows from $(*2),(*4)$ and Corollary \ref{corl-5}.
	\end{proof}
	
	\begin{proof}[The proof of Lemma \ref{lm-8}]
		For each $x \in S\left(A^{+}, A^{-}, b\right)$, Lemma \ref{lm-5} implies that $x_{j_{0}} \in I_{j_{0}}=\{k\}$, i.e., $x_{j_{0}}=k$. Now, from Lemma \ref{lm-6}, if $x \in S_{\left\{i_{0}\right\}}\left(A^{+}, A^{-}, b\right)$, then we have $x_{j} \in I_{i_{0} j}$, $\forall j \in \mathscr{J}\left(*1\right)$. Particularly, $x_{j_{0}} \in I_{i_{0} j_{0}} \subseteq I_{j_{0}}=\{k\}$. On the other hand, since $k \in S_{i_{0}{j_{0}}}^{\prime}$, then we have $S_{i_{0} j_{0}}^{\prime} \neq \varnothing$, that together with $S_{i_{0} j_{0}}^{\prime}=S_{i_{0} j_{0}} \cap I_{j_{0}}$ and $I_{j_{0}}=\{k\}$ imply $S_{i_{0} j_{0}}^{\prime}=\{k\}$. Consequently, $x_{j_{0}} \in S_{i_{0} j_{0}}^{\prime}\left(*2\right)$. Now, $\left(*1\right)$,$\left(*2\right)$ and Corollary \ref{corl-5} imply that $x$ satisfies the $i_{0}$'th equation.
	\end{proof}
	
	\begin{proof}[The proof of Lemma \ref{lm-9}]
		From $x \in S_{\left\{i_{0}\right\}}\left(A^{+}, A^{-}, b\right)$ and Lemma \ref{lm-6}, we conclude that $x_{j} \in I_{i_{0} j}$, $\forall j \in \mathscr{J}\left(*1\right)$. Also, from Lemma \ref{lm-5}, there exists at least one $j_{i} \in \mathscr{J}$ such that $x_{j_{i}} \in S_{ij_{i}}^{\prime}$. But, since $S_{i j}^{\prime} \subseteq S_{i_{0} j}^{\prime}$, it follows that $x_{j_{i}} \in S_{i_{0} j_{i}}^{\prime}(*2)$. So, $(*1),(*2)$ and Corollary \ref{corl-5} imply that $x$ satisfies the $i_{0}$'th equation.
	\end{proof}
	
	\begin{proof}[The proof of Lemma \ref{lm-10}]
		If $x \in S\left(A^{+}, A^{-}, b\right)$, then from Lemma \ref{lm-5}, there exists at least one $j \in \mathscr{J}$ such that $x_{j} \in S_{i_{0} j}^{\prime}$. But, since $\mathscr{J}_{i_{0}}=\left\{j_{0}\right\}$, we have necessarily $x_{j_{0}} \in S_{i_{0} j_{0}}^{\prime}=\{k\}$, that means $x_{j_{0}}=k$. Moreover, if $x \in S_{\{i\}}\left(A^{+}, A^{-}, b\right) \cap\left\{x \in[0,1]^{n}: x_{j_{0}}=k\right\}$, then $x_{j_{0}}=k \in S_{i j_{0}}^{\prime}$ and also from Lemma \ref{lm-6}, $x_{j} \in I_{i j}$, $\forall j \in \mathscr{J}$. So, Corollary \ref{corl-5} implies that $x$ satisfies the $i$'th equation.
	\end{proof}
	
	\begin{proof}[The proof of Lemma \ref{lm-11}]
		Assume that $x(e) \in S_{\left\{i_{0}\right\}}\left(A^{+}, A^{-}, b\right)$. To prove the lemma, it is sufficient to show that $x(e)$ also satisfies the $i_{0}$'th equation. First, from Lemma \ref{lm-6}, it follows that $x(e)_{j} \in I_{i_{0} j}$, $\forall j \in \mathscr{J}\left(*1\right)$. On the other hand, from $\left|S_{i_{0}j_{0}}^{\prime}\right|=2$ and Corollary \ref{corl-4}, we have $S_{i_{0} j_{0}}^{\prime}=\left\{L_{j_{0}}, U_{j_{0}}\right\}$. Moreover, based on (\ref{eq-12}), it follows that either $x(e)_{j_{0}}=L_{j_{0}}$ (if $\mathscr{I}_{j_{0}}(e)=\varnothing$) or $x(e)_{j_{0}}=\min \left\{\bigcap_{i \in \mathscr{I}_{j_{0}}(e)} S_{ij_{0}}^{\prime}\right\}$ (if $\mathscr{I}_{j_{0}}(e) \neq \varnothing$ ). In the former case, it is clear that $x(e)_{j_{0}}=L_{j_{0}} \in S_{i_{0} j_{0}}$. Otherwise, in the latter case, it follows that $\bigcap_{i \in \mathscr{I}_{j_{0}}(e)} S_{ij_{0}}^{\prime} \neq \varnothing$ (Corollary \ref{corl-6}) and $\min \left\{\bigcap_{i \in \mathscr{I}_{j_{0}}(e)} S_{ij_{0}}^{\prime}\right\} \in\left\{L_{j_{0}}, U_{j_{0}}\right\}$ (Corollary \ref{corl-4}), and therefore in any case we have $x(e)_{j_{0}} \in S_{i_{0} j_{0}}\left(*2\right)$. As a result, $\left(*1\right)$, $\left(*2\right)$ and Corollary \ref{corl-5} imply that $x(e)$ also satisfies the $i_{0}$'th equation.
	\end{proof}
	
	\begin{proof}[The proof of Lemma \ref{lm-12}]
		Let $e \in E$. So, from (\ref{eq-12}), we have either $x(e)_{j_{0}}=L_{j_{0}}$ or $x(e)_{j_{0}}=\min \left\{\bigcap_{i \in \mathscr{I}_{j_{0}}(e)} S_{ij_{0}}^{\prime}\right\}$. In the latter case, it is concluded that $\bigcap_{i \in \mathscr{I}_{j_{0}}(e)} S_{\mathrm{ij}_{0}}^{\prime} \neq \varnothing$ (Corollary \ref{corl-6}), which requires $\mathscr{I}_{j_{0}}(e) \subseteq \mathscr{I}_{j_{0}}$. Hence, by our assumption, we have $L_{j_{0}} \in S_{ij_{0}}^{\prime}$, $\forall i \in \mathscr{I}_{j_{0}}(e)$, that implies $\min \left\{\bigcap_{i \in \mathscr{I}_{j_{0}}(e)} S_{ij_{0}}^{\prime}\right\}=L_{j_{0}}$. So, in any case, $x(e)_{j_{0}}=L_{j_{0}}$. Now, let $L_{j_{0}} \in S_{i_{0} j_{0}}^{\prime}$ and $x(e) \in S_{\left\{i_{0}\right\}}\left(A^{+}, A^{-}, b\right)$. So, from Lemma \ref{lm-6}, it follows that $x(e)_{j} \in I_{i_{0} j}$, $\forall j \in \mathscr{J} ~(*)$. Moreover, we have either $x(e)_{j_{0}}=L_{j_{0}}$ or $x(e)_{j_{0}}=\min \left\{\bigcap_{i \in \mathscr{I}_{j_{0}}(e)} S_{ij_{0}}^{\prime}\right\}=L_{j_{0}}$, where the latter equality results from the assumption that $L_{j_{0}} \in S_{ij_{0}}^{\prime}$, $\forall i \in \mathscr{I}_{j_{0}}$. Hence, $x(e)_{j_{0}} \in S_{i_{0} j_{0}}^{\prime}$, which together with $(*)$ and Corollary \ref{corl-5} imply that $x(e)$ satisfies the $i_{0}$'th equation.
	\end{proof}
	
	\begin{proof}[The proof of Lemma \ref{lm-13}]
		\textbf{(a)} Based on the definition of $F_{0}$, it is clear that $F_{0} \subseteq F$. To prove the inclusion $S^{*} \subseteq F_{0}$, we shall show that for each $e \in E-E_{0}$, there exists some $e_{0} \in E_{0}$ such that $\sum_{j \in \mathscr{J}} c_{j} x\left(e_{0}\right)_{j} \leq \sum_{j \in \mathscr{J}} c_{j} x(e)_{j}$. For this purpose, let $e \in E-E_{0}$. So, $\mathscr{I}_{j_{1}}(e) \neq \varnothing$. Now, consider the function $e_{0}$ as follows:
		\begin{equation}
			e_{0}(i)= \begin{cases}j_{2} & , i \in \mathscr{I}_{j_{1}}(e)  \\ e(i) & , i \notin \mathscr{I}_{j_{1}}(e)\end{cases}
		\end{equation}
		Hence, $\mathscr{I}_{j_{1}}\left(e_{0}\right)=\varnothing$ that implies $e \in E_{0}$, $\mathscr{I}_{j_{1}}\left(e_{0}\right) \subseteq \mathscr{I}_{j_{1}}(e)$ and $x\left(e_{0}\right)_{j_{1}}=L_{j_{1}}$ (from (\ref{eq-12})). However, $\bigcap_{i \in \mathscr{I}_{j_{1}}(e)} S_{i j_{1}}^{\prime} \neq \varnothing$ (Corollary \ref{corl-6}) and $\min \left\{\bigcap_{i \in \mathscr{I}_{j_{1}}(e)} S_{ij_{1}}^{\prime}\right\} \in\left\{L_{j_{1}}, U_{j_{1}}\right\}$ (Corollary \ref{corl-4}), which requires $x(e)_{j_{1}}=v \in\left\{L_{j_{1}}, U_{j_{1}}\right\}$. So, $c_{j_{1}} x\left(e_{0}\right)_{j_{1}}=c_{j_{1}} L_{j_{1}} \leq c_{j_{1}} v=c_{j_{1}} x(e)_{j_{1}}\left(*1\right)$. Moreover, based on (16), we have $\mathscr{I}_{j}\left(e_{0}\right)=\mathscr{I}_{j}(e)$, $\forall j \in \mathscr{J}-\left\{j_{1}, j_{2}\right\}$, which means $x\left(e_{0}\right)_{j}=x(e)_{j}$, and therefore $c_{j} x\left(e_{0}\right)_{j}=c_{j} x(e)_{j}$, $\forall j \in \mathscr{J}-\left\{j_{1}, j_{2}\right\} (*2)$. Subsequently, the assumption $\bigcap_{i \in \mathscr{I}_{j_{2}}} S_{ij_{2}}^{\prime}=\left\{L_{j_{2}}\right\}$ together with Remark \ref{rmk-6} imply that $x\left(e_{0}\right)_{j_{2}}=\min \left\{\bigcap_{i \in \mathscr{I}_{j_{2}}\left(e_{0}\right)} S_{ij_{2}}^{\prime}\right\}=L_{j_{2}}$. Also, if $\mathscr{I}_{j_{2}}(e)=\varnothing$ we have $x(e)_{j_{2}}=L_{j_{2}}$ (from (\ref{eq-12})) and if $\mathscr{I}_{j_{2}}(e) \neq \varnothing$ we have $x(e)_{j_{2}}=\min \left\{\bigcap_{i \in \mathscr{I}_{j_{2}}(e)} S_{ij_{2}}^{\prime}\right\}=L_{j_{2}}$ (Remark \ref{rmk-6} and (\ref{eq-12})). Hence, $x\left(e_{0}\right)_{j_{2}}=x(e)_{j_{2}}$ that requires $c_{j_{2}} x\left(e_{0}\right)_{j_{2}}=c_{j_{2}} x(e)_{j_{2}} \left(*3\right)$. Now, from $\left(*1\right)$, $\left(*2\right)$ and $(*3)$, $c_{j_{1}} x\left(e_{0}\right)_{j_{1}}+c_{j_{2}} x\left(e_{0}\right)_{j_{2}}+\sum_{j \in \mathscr{J}-\left\{j_{1}, j_{2}\right\}} c_{j} x\left(e_{0}\right)_{j} \leq c_{j_{1}} x(e)_{j_{1}}+c_{j_{2}} x(e)_{j_{2}}+\sum_{j \in \mathscr{J}-\left\{j_{1}, j_{2}\right\}} c_{j} x(e)_{j}$.
		\textbf{(b)} It is quite similar to the proof of Part (a).
	\end{proof}
	
	\section*{Acknowledgement}
	The authors wish to express their appreciation for several excellent suggestions
	for improvements in this paper made by the referees.\\

\bibliography{refrences}

\end{document}